\documentclass[12pt]{article}
\usepackage{amsmath,amssymb,amscd,latexsym,euscript}
\usepackage[matrix,arrow,curve]{xy}
\newcommand{\coLim}{\underrightarrow{\lim}}
\newcommand{\Lim}{\underleftarrow{\lim}}
\newcommand{\hd}{{\rm hd\,}}
\newcommand{\Ld}{{\rm Ld\,}}
\newcommand{\emb}{{\rm emb\,}}
\newcommand{\Hom}{{\rm Hom}}
\newcommand{\Lan}{{\rm Lan}}
\newcommand{\Add}{{\rm Add\,}}
\newcommand{\Set}{{\rm Set}}
\newcommand{\Ab}{{\rm Ab}}
\newcommand{\Ob}{{\rm Ob\,}}
\newcommand{\Mor}{{\rm Mor\,}}
\newcommand{\init}{{\rm init\,}}

\newcommand{\Imm}{{\rm Im\,}}
\newcommand{\Ker}{{\rm Ker\,}}

\newcommand{\Mat}{{\rm Mat\,}}
\newcommand{\Tor}{{\rm Tor\,}}

\newcommand{\NN}{{\,\mathbb N}}

\newcommand{\RR}{{\,\mathbb R}}

\newcommand{\mC}{{\cal \, C}}
\newcommand{\mD}{{\cal \, D}}
\newcommand{\mK}{{\cal \, K}}
\newcommand{\mHS}{{\cal \, HS}}

\newcommand{\ZZ}{{\,\mathbb Z}}
\newcommand{\II}{{\,\mathbb I}}
\newcommand{\mA}{{\mathcal A}}
\newcommand{\mB}{{\mathcal B}}

\newcommand{\fA}{{\mathfrak{A}}}

\newcommand{\fS}{{\mathfrak{S}}}
\newcommand{\fF}{{\mathfrak{F}}}

\newtheorem{theorem}{\bf Theorem}[section]
\newtheorem{lemma}[theorem]{\bf Lemma}
\newtheorem{proposition}[theorem]{\bf Proposition}
\newtheorem{corollary}[theorem]{\bf Corollary}
\newtheorem{definition}{\sc Definition}[section]
\newtheorem{example}[definition]{\sc Example}
\newtheorem{remark}[definition]{\sc Remark}

\def\leq{\leqslant}
\def\geq{\geqslant}

\title
{
On the Cubical Homology Groups of Free Partially Commutative Monoids 
}
\author{Ahmet A. Husainov}
\date{}

\begin{document}

\maketitle

\begin{abstract}
We study a Leech homology of a 
locally bounded free partially commutative monoid $M(E,I)$.
Given 
a contravariant natural system of 
abelian groups $F$ on  $M(E,I)$
we build a precubical set 
$T(E,I)$ with a homological system of abelian groups  $\overline{F}$ 
  and prove that the Leech homology groups $H_n(M(E,I),F)$ are isomorphic
  to the
cubical homology groups $H_n(T(E,I),\overline{F})$,
 $n\geq 0$.
As a consequence we have confirmed a conjecture that 
if the free partially commutative 
monoid does not contain $>n$ mutually commuting generators, then its the 
homological dimension $\leq n$.
We have built the complexes of finite length 
for a computation of the Leech homology of  
such monoids and the Hochschild homology of their monoid rings.
The results are applied to the homology of asynchronous 
transition systems. We give the positive answer to a question that 
the homological dimension of the asynchronous system does not greater 
than the maximal number of its  
mutually independent events.
We have built the complex for computing the integral homology groups 
of an asynchronous transition system by the Smith normal form of integer 
matrices.
\end{abstract}

2000 Mathematics Subject Classification 16E40, 18B20, 18G10, 18G20,
 18G35, 18G40, 20M50, 55U99, 68Q85


\tableofcontents

\section*{Introduction}

This work contains three sections. 
In the first section we study  homology groups 
of precubical sets with coefficients 
in homological systems of abelian groups.
It is known \cite[Application 2, Prop. 4.2]{gab1967}
that the homology groups $H_n(X,F)$ 
of any simplicial set $X$
with coefficients in a homological system of abelian groups $F$ 
are isomorphic to the values $\coLim_n^{(\Delta/X)^{op}}F$ 
 of the left satellites of the colimit where 
$\Delta/X$ is the category of singular simplexes of $X$.
Similar assertion is true for the homological systems 
on semisimplicial sets 
(see  \cite[Prop. 1.4]{X1997} for the dual assertion). 
It does not true in case of cubical sets (see Example \ref{concub}).
Nevertheless the main result of the first section 
(Theorem \ref{comcub}) shows that it is true for homological systems 
on the precubical sets. This allows us to use 
the homology theory of small categories.
Using Oberst's Theorem \cite[Theorem 2.3]{obe1968}, 
we find the criterion  
under which the homomorphisms of homology groups of precubical sets 
corresponding to the precubical morphism are isomorphisms.
We will build spectral sequences 
of a locally ordered covering 
(Cor. \ref{locdir}) and of a morphism of precubical sets
(Cor. \ref{spmor}).

Second section is devoted to the homological dimension of monoids
satisfying the following 
\begin{definition}
Let $E$ be a set and $I\subseteq E\times E$ an irreflexive and 
symmetric relation. 
A monoid given by a set of generators  $E$ and relations $ab=ba$   
 for all  $(a,b)\in I$ is called {\em free partially commutative} and  
denoted by $M(E,I)$.
If  $(a,b)\in I$ then the members $a,b\in E$ are said to be 
{\em commuting generators}.
\end{definition}
This definition is varied with ordinary \cite{die1997}. We do not demand 
that  the set $E$ is finite. 
It was shown in  \cite{X20042} that if the monoid $M(E,I)$ does 
not contain triples of mutually commuting generators, then its 
homological dimension $\leq 2$. 
We put forward a conjecture in \cite{X2005} that if 
$M(E,I)$ does not contain $n$--tuples of mutually commuting generators,
then its
homological dimension  $\leq n$. 
This conjecture was confirmed in the case of finite $E$ 
by L. Polyakova \cite{pol2006}. We have confirmed it 
in the general case (Cor. \ref{hyptes1}).

We will prove that Leech homology 
groups of free partially commutative 
monoids are isomorphic to the cubical homology groups.
This result we use to the building complexes for computing
the homology groups
of monoids and Hochschild homology of monoid rings.

Third section is devoted to the homology of $M(E,I)$--sets.
M. Bednarczyk \cite{bed1988} has introduced
{\em asynchronous transition systems}
to the modelling the concurrent processes. 
In \cite{X20032} 
it was proved that the category of asynchronous
transition systems admits an inclusion into the category of
pointed sets over free partially commutative monoids.
Thus asynchronous
transition systems may be considered as 
$M(E,I)$--sets. It allows us to build the 
functors from the category of asynchronous transition systems 
into the category of the abelian group and homomorphisms.
Their values are homology groups of asynchronous transition systems.
It was shown in \cite{X20032} (see \cite{X20042} as well) 
that these groups are isomorphic to the homology groups of 
$M(E,I)$ with coefficients in some right 
$M(E,I)$--modules. 
Goubault \cite{gou1995} and Gaucher \cite{gau1999}, \cite{gau2000}
have given a definition of homology groups for 
higher dimensional automata which are precubical sets with 
additional structure. 
There is a functor from the 
category of asynchronous transition systems into the category 
of precubical sets. 
In the third section we prove that homology groups  
of asynchronous transition systems, defined in \cite{X20032},
are isomorphic to the homology groups of the corresponding precubical sets.
In the work \cite{X20042} it was proved that if
an asynchronous transition system
  does not contain triples of independent events, then its homological 
dimension $\leq 2$. It was given a conjecture that 
if the asynchronous transition system  
  does not contain $n$--tuples of independent events, then its homological 
dimension $\leq n$ \cite[Open Problem 2]{X20042}. 
We have confirmed this conjecture.

We use the following notations:

\noindent
\medskip
\begin{tabular}{p{2.5cm} p{10cm}}
$\Set$ & is the category of sets and maps.\\ 
$\Set_*$ & is the category of pointed sets: each its object is a set
$X$ with
a selected element, written  $*$ and called the ``base point" ; 
its morphisms are maps $X\rightarrow Y$ 
which send the base point of $X$
to that of $Y$. Such maps are called {\em based}.
Following \cite{win1995} we identify
the category of sets and partial functions
with $Set_*$.  \\
$\Ab$ & is the category of abelian groups and homomorphisms.\\
$L: \Set \rightarrow \Ab$ & is the left adjoint to the forgetful functor 
$U: Ab\rightarrow Set$; $L$ assign to every set $E$ the free 
abelian group $L(E)$ generated by $E$ and to a map 
$f: E_1\rightarrow E_2$ the canonical   
homomorphism $L(f): L(E_1)\rightarrow L(E_2)$ such that 
$L(f)(e)=f(e)$ for all $e\in E_1$.\\
$\II$ & is the totally ordered set $\{0,1\}$ with the order relation
 $0\leq 1$.\\
$\ZZ$ & is the additive group of integers.\\
$\NN$ & is the set of non-negative integers or the additive monoid 
\{ $0$, $1$, 
$2$, $3$, $\cdots$ \} or the free multiplicative monoid  
$\{1, a, a^2, a^3, \cdots \}$.\\
$\RR$ & is the set of reals.
\end{tabular}

For any category $\mA$, denote by $\mA^{op}$ the dual category. 
For any pair  $a,b\in Ob\mA$, denote by 
$\mA(a,b)$ the set of all morphisms  $a\rightarrow b$.  
Given  a small category $\mC$ we denote by $\mA^{\mC}$ the category of 
functors $\mC\rightarrow \mA$ and natural transformations.

We will be consider any monoid  as the small category 
with one object.
This exert influence on our terminology. In particular a right $M$--set $X$
will be considered and denoted as a functor 
 $X: M^{op}\rightarrow \Set$ (the value of $X$ at the unique object  
will be denoted by $X(M)$ or shortly $X$). 
Morphisms of right $M$--sets are natural transformations.

\section{Homology of categories and precubical sets}

\subsection{Precubical sets}

A {\em precubical set} 
$X= (X_n, \partial_i^{n,\varepsilon})$ 
is a sequence of sets
$(X_n)_{n\in \NN}$ with a family of maps 
$\partial_i^{n,\varepsilon}: X_n \rightarrow X_{n-1}$, 
defined for  $1\leq i\leq n$, $\varepsilon\in \{0,1\}$, for which 
the following diagrams is commutative for all
$\alpha,\beta \in \{0,1\}$, $n\geq 2$ и $1\leq i< j\leq n$:
$$
\begin{CD}
X_n @>\partial_j^{n,\beta}>> X_{n-1}\\
@V{\partial_i^{n,\alpha}}VV @VV\partial_i^{n-1,\alpha}V\\
X_{n-1} @>>\partial_{j-1}^{n-1,\beta}> X_{n-2}
\end{CD}
$$
For example, if $X_n=\emptyset$ and  $n>1$, then the maps 
$\partial_1^{1,0}: X_1\rightarrow X_0$ and 
$\partial_1^{1,1}: X_1\rightarrow X_0$ will be determine 
a directed graph.

\medskip
\noindent
{\bf Precubical sets as functors.}
Let 
$\Box_+$ be a category consisting of the finite  sets
 $\II^n=\{0,1\}^n$ ordered as the Cartesian power of $\II$.
Any morphism of $\Box_+$ is defined as an ascending map which admits 
a decomposition  of the form
$\delta_i^{k,\varepsilon}: \II^{k-1}\rightarrow \II^k$
where
$$
\delta_i^{k,\varepsilon}(x_1, \cdots, x_{k-1})=
(x_1, \cdots, x_{i-1}, \varepsilon, x_i, \cdots, x_{k-1}), \quad
 \varepsilon\in \II, 
1 \leq i \leq k.
$$
Every morphism  $f\in \Box_+(\II^m, \II^n)$ 
has a decomposition 
$f=\delta_{j_{n-m}}^{n,\varepsilon_{n-m}}\cdots
\delta_{j_1}^{m+1,\varepsilon_1}$,
$1\leq j_1 < \cdots < j_{n-m}\leq n$
because of the equations
$\delta_j^{n,\beta}\delta_i^{n-1,\alpha}= 
\delta_i^{n,\alpha}\delta_{j-1}^{n-1,\beta}$
for $1\leq i< j \leq n$, $\alpha\in \{0,1\}$, $\beta\in \{0,1\}$.
This allows us to define a precubical set as a functor
$X: \Box_{+}^{op} {\rightarrow}\Set$ with values 
$X(\II^n)=X_n$ at $\II^n\in Ob(\Box_{+})$
and  
$X(f)= 
\partial_{j_1}^{m+1,\varepsilon_1}\dots
\partial_{j_{n-m}}^{n,\varepsilon_{n-m}}$
at  
$f=\delta_{j_{n-m}}^{n,\varepsilon_{n-m}}\dots
\delta_{j_1}^{m+1,\varepsilon_1}: \II^m\rightarrow \II^n$.
A morphism of precubical sets is defined as a natural transformation.
So we get the category of precubical sets $\Set^{\Box_+^{op}}$.

\medskip
\noindent
{\bf Integral homology groups of precubical sets.}
\begin{definition}\label{homolcub}
Let 
$X= (X_n, \partial_i^{n,\varepsilon})$ 
be a precubical set and $C_*(X)$ the chain complex of abelian groups 
$C_n(X)= L(X_n)$, $n\geq 0$, with $C_n(X)= 0$ for $n<0$.
Its differentials $d_n: C_n(X)\rightarrow C_{n-1}(X)$ are defined by
$$
	d_n = \sum_{i=1}^n(-1)^i(L(\partial_i^{n,1})-L(\partial_i^{n,0})).
$$
A {\em $n$-th integral  
homology group of the precubical set  X}
is the group  $H_n(C_*(X))$.
\end{definition}

\medskip
\noindent
{\bf Cubical subsets of Euclidian space.}
Consider precubical sets which integral homology groups are studied in
 \cite{kac2000}-\cite{kac2003}.

For each integer $l\in\ZZ$ the closed intervals 
$[l,l+1]\subset \RR$ and $\{l\}=[l,l]\subset \RR$ are called  
{\em elementary intervals}. The elementary intervals $[l,l+1]$
are  {\em nondegenerate} and $[l,l]$ {\em degenerate}.
\begin{definition}
An {\em elementary cube}
is a Cartesian product of elementary intervals
$Q= I_1\times I_2 \times \dots \times I_n \subset \RR^n$.
The number of nondegenerate intervals among  $I_1, I_2, \dots, I_n$ is 
denoted by $\dim Q$.
Denote the dimension of Euclidian space  by $\emb Q= n$.
\end{definition}
\begin{definition}
A {\em  cubical subset of Euclidian space} 
$\RR^n$ is the  union of arbitrary elementary cubes $Q_1, Q_2, \cdots$, 
for which $\emb Q_i = n$.
\end{definition}
For any cubical subset $X\subseteq \RR^n$ 
and $m\geq 0$ denote by 
$\mK_m(X)$ the set of elementary cubes $Q\subseteq X$ such that $\dim Q= m$.

\begin{example}
$X= [0,1]\times[0,2] \cup [1,2]\times [0,0] \cup [2,2]\times [0,1] 
\cup [1,2]\times [1,1] $ (Pic. \ref{semicub} ). 

$\mK_0(X) = \{ [0,0]\times [0,0], [1,1]\times [0,0], [2,2]\times [0,0],$
		$[0,0]\times[1,1], [1,1]\times[1,1], [2,2]\times[1,1],$
	$[0,0]\times[2,2], [1,1]\times[2,2]    \}$.

$\mK_1(X) = \{ [0,1]\times [0,0], [1,2]\times [0,0], [0,1]\times [1,1],$
		$[1,2]\times[1,1], [0,1]\times[2,2], [0,0]\times[0,1],$
  $[1,1]\times[0,1], [2,2]\times[0,1], [0,0]\times [1,2], [1,1]\times [1,2]$.

$\mK_2(X) = \{ [0,1]\times [0,1], [0,1]\times [1,2]$.

\begin{figure}[ht]
\begin{center}
\begin{picture}(240,160)

\put(20,40){\vector(1,0){180}} \put(202,32){$x_1$}
\put(40,20){\vector(0,1){120}} \put(25,140){$x_2$}

\multiput(40,41)(0,5){17}{\line(1,0){40}}
\put(40,80){\line(1,0){80}}
\put(40,120){\line(1,0){40}}
\put(80,40){\line(0,1){80}}
\put(120,40){\line(0,1){40}}
\put(30,30){$0$}
\put(80,30){$1$}
\put(120,30){$2$}
\put(30,80){$1$}
\put(30,120){$2$}
\end{picture}
\end{center}
\caption{\label{semicub} A cubical subset of the space $\RR^2$.}
\end{figure}
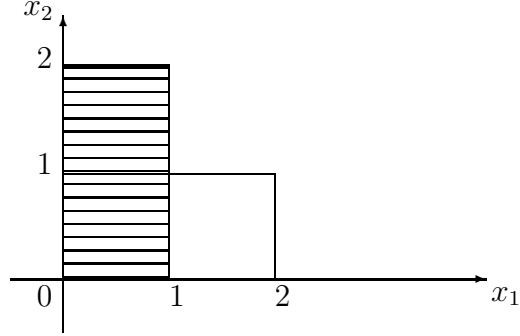

\end{example}

Let  $X\subseteq \RR^n$ be a cubical subset.
For any elementary cube  
$Q=I_1\times I_2\times \cdots \times I_n \in \mK_m(X)$, 
where $I_{k_j}= [a_j,b_j]$ is a $j$-th nondegenerate elementary 
interval, we let
$$
\partial_j^{m,0} Q = I_1\times I_2\times \cdots
\times I_{k_j-1}\times \{a_j\}
\times I_{k_j+1}\times \cdots \times I_n\quad,
$$
$$
\partial_j^{m,1} Q = I_1\times I_2\times \cdots
\times I_{k_j-1}\times \{b_j\}
\times I_{k_j+1}\times \cdots \times I_n\quad, 
$$
where $1\leq j\leq \dim Q =m$.
We have got the precubical set
$\mK(X)=(\mK_m(X), \partial_j^{m,\varepsilon})$.
Consequently, the cubical subsets of Euclidian spaces may be 
regarded as precubical sets.
Easy to see (from \cite[Corollary 2.3]{kac2003}, for example), 
that the homology groups of precubical sets corresponding 
to the cubical subsets of Euclidian spaces are isomorphic 
to the cubical homology groups defined in  \cite{kac2003}.
In particular it is true the following assertion:

\begin{lemma}\label{KMM}
Let $C_*(\mK([0,1]^n))$ be a chain complex of abelian 
groups
$C_m(\mK([0,1]^n))=L(\mK_m([0,1]^n))$ generated by cubes
$Q\subseteq [0,1]^n$ for which $\dim Q=m$. 
Let $Q= I_1\times \cdots \times I_n$ be an elementary cube which equals 
a product of $m$ elementary nondegenerate intervals 
 $I_{k_1}=\cdots = I_{k_m}=[0,1]$, $1\leq k_1 < k_2 < \cdots < k_m \leq n$,
and $n-m$ degenerate elementary intervals. Define the differentials 
$d_m$  of complex $C_*(\mK([0,1]^n))$ 
as
$$
d_m Q = \sum\limits_{j=1}^m 
(-1)^j(L(\partial_{j}^{m,1}) Q - L(\partial_{j}^{m,0}) Q).
$$
Then 
$$
H_q(C_*(\mK([0,1]^n)))=
\left\{
\begin{array}{ccc}
\ZZ, & $if$ & q=0,\\
0, & $if$ & q>0.\\
\end{array}
\right.
$$
\end{lemma}

This assertion is proved in \cite[Prop. 4.61]{kac2000} for the most 
general case where an arbitrary elementary
cube of Euclidian space is taken
 instead of  $[0,1]^n$.

\subsection{Smith normal form and calculating the homology groups}

Let
$G'\stackrel\alpha\rightarrow G\stackrel\beta\rightarrow G''$ be
homomorphisms of abelian groups such that $\beta\circ\alpha=0$.
Then we can consider the homology group $\Ker(\beta)/\Imm(\alpha)$.
Let us describe a method of calculation for this homology group when
$G'$, $G$, and $G''$ are finitely generated free abelian groups.
\begin{proposition}{\em \cite[Theorem III.4(43)]{nod1999}}
Let
$$
A =
\begin{pmatrix}
a_{11}& a_{12}& \cdots & a_{1n}\\
a_{21}& a_{22}& \cdots & a_{2n}\\
\vdots& \vdots & \ddots & \vdots\\
a_{m1}& a_{m2}& \cdots & a_{mn}
\end{pmatrix}
$$
be a matrix with integer entries $a_{ij}\in\ZZ$. Then there
{is} an $m\times{m}$ matrix $T$ and {an}   $n\times n$
matrix $S$ with integer entries such that:

(i) $det(T)=\pm{1}$, $det(S)=\pm{1}$;

(ii) $A=T\circ D(A)\circ S$ for a natural number
$k\geq 0$ and  $m\times n$ matrix
$$
D(A) =
\begin{pmatrix}
d_{1}& 0 & \cdots & 0 &0& \cdots & 0\\
0 & d_{2}& \cdots & 0 &0& \cdots & 0\\
\vdots& \vdots & \ddots &\vdots& \vdots & & \vdots\\
0 & 0 & \cdots & d_k & 0& \cdots & 0\\
0 & 0 & \cdots & 0 & 0&\cdots  & 0\\
\vdots& \vdots &  &\vdots& \vdots & & \vdots\\
0 & 0 & \cdots & 0 &0& \cdots & 0
\end{pmatrix}
$$
all entries of which are equal to $0$ except the diagonal numbers
$d_1\leq d_2\leq \cdots \leq d_k$ which satisfy
that $d_i$ divides $d_{i+1}$ for all $1\leq{i}\leq{k-1}$.
\end{proposition}
The matrix  $D(A)$ is said to be a {\em Smith normal form of}
$A$.
This form is used for the computation of the homology
groups of simplicial complexes in \cite{sei2001}. We refer the
reader to \cite{hav1997}, where  an algorithm is presented
for computing the Smith normal form of an integer matrix, which
performs well in practice.
 There are packages such as {\sc GAP} for the computation of the
Smith normal form.

\begin{proposition}\label{Smith}
Let a homomorphism $\alpha:\ZZ^n \rightarrow \ZZ^m$ be given by a
$m\times{n}$ matrix $A$ and  $\beta:\ZZ^m \rightarrow \ZZ^p$ by
$p\times{m}$ matrix $B$. Suppose the sequence
$\ZZ^n\stackrel\alpha\rightarrow \ZZ^m\stackrel\beta\rightarrow \ZZ^p$
satisfies to
$\beta\circ\alpha=0$. If
$d_1$, $d_2$, $\cdots$, $d_k$ are the nonzero diagonal entries
of the Smith normal form
 $D(A)$, then
\begin{equation}\label{iso}
\Ker\beta/\Imm\alpha \cong \ZZ/d_1\ZZ \times
\ZZ/d_2\ZZ \times \cdots \times \ZZ/d_k\ZZ
\times \ZZ^{m-k-b}
\end{equation}
where  $b$ is the rank of the matrix $B$ and may be
computed as the number of nonzero entries of
$D(B)$.
\end{proposition}

\subsection{Homology of small categories}

In this section it is considered 
homology groups of small categories 
  $\mC$
with coefficients in functors $F: \mC \rightarrow \Ab$.

\medskip
\noindent
{\bf Homology of categories 
and derived functors of the colimit.}
\begin{definition}\label{defhomolcat}
{\em Homology groups of a small category $\mC$
with coefficients in a functor
 $F:\mC\rightarrow \Ab$} are homology groups of the chain complex
 $C_*(\mC,F)$
consisting of the abelian groups
\begin{displaymath}
 C_n({\mC},F) = \bigoplus_{c_0 \rightarrow \cdots \rightarrow c_n}
F(c_0), \quad n \geq 0,
\end{displaymath}
with differentials
$d_n= \sum\limits_{i=0}^{n}(-1)^i d^n_i: 
C_n(\mC,F) \rightarrow C_{n-1}(\mC,F)$, $n>0$,
where  for each 
$$
(c_0 \stackrel{\alpha_1}\rightarrow c_1 \stackrel{\alpha_2}\rightarrow
\cdots \stackrel{\alpha_{n}}\rightarrow c_{n}, a)
 \in 
\bigoplus_{c_0 \rightarrow \cdots \rightarrow c_n}F(c_0), \quad a\in F(c_0)
$$
$d^n_i$ are defined by
$$
d^n_i(c_0 \rightarrow \cdots \rightarrow c_n, a) =
\left\{
\begin{array}{ll}
(c_0 \rightarrow \cdots \rightarrow \widehat{c_i} \rightarrow 
\cdots \rightarrow c_n, a) \quad , & \mbox{for} ~ 1 \leq i \leq n\\
(c_1 \rightarrow \cdots \rightarrow c_n,
F(c_0\stackrel{\alpha_1}\rightarrow c_1)(a) ) ~~,
& \mbox{for}~ i = 0.
\end{array}
\right.
$$
\end{definition}
It is well known \cite[Prop. 3.3, Application 2]{gab1967} that 
there exists an isomorphism of left satellites of the colimit  
$\coLim^\mC: \Ab^\mC \rightarrow \Ab$ 
and the functors  
 $H_n(C_*(\mC,-)): \Ab^\mC \rightarrow \Ab$.
Since the category $\Ab^\mC$ has enough projectives,
these satellites are natural isomorphic to the left derived functor of 
$\coLim^{\mC}: \Ab^{\mC}\rightarrow \Ab$.
 Denote the values $H_n(C_*(\mC,F))$ of the satellites at $F\in \Ab^\mC$
   by $\coLim_n^\mC F$.

Let  $\mC$ be a small category. Denote by
$\Delta_\mC\ZZ: \mC\rightarrow \Ab$, or shortly $\Delta\ZZ$,  the functor 
which has constant values  $\Delta\ZZ(c)= \ZZ$ at $c\in Ob(\mC)$ 
and $\Delta\ZZ(\alpha)=1_\ZZ$ at $\alpha \in Mor(\mC)$. 
{\em Integral homology groups} of $\mC$
are the 
groups $\coLim_n^\mC\Delta_\mC\ZZ$.
\begin{example}
Let $pt=[0]$ be a category consisting of one object $0$
and one morphism $1_0$.
Then for every $n\geq 0$ the abelian group $C_n(pt,\Delta_{pt}\ZZ)$
is generated by the element $e_n$ which is equal to the sequence 
of  $n$ identity morphisms 
$0\stackrel{1_0}\rightarrow 0 \stackrel{1_0}\rightarrow 
\cdots \stackrel{1_0}\rightarrow 0$.
Since $d_i^n(e_n)=e_{n-1}$,  
$$
d_n(e_n)=
\left\{
\begin{array}{ccc}
e_{n-1}, & $if$ & n $ is even$,\\
0, & $if$ & n $ is odd$.\\
\end{array}
\right.
$$
Consequently the chain $C_*(\mC,\Delta\ZZ)$ is isomorphic to the chain
$$
0 \leftarrow \ZZ \stackrel{0}\leftarrow \ZZ \stackrel{1}\leftarrow 
\ZZ \stackrel{0}\leftarrow \ZZ \stackrel{1}\leftarrow  \cdots ,
$$
which $n$-th homology groups are equal to  $0$ if $n>0$, and 
$H_0(C_*(pt,\Delta\ZZ))\cong \ZZ$.
It follows that
$$
\coLim_n^{pt}\Delta\ZZ=
\left\{
\begin{array}{ccc}
\ZZ, & $if$ & n=0,\\
0, & $if$ & n>0.\\
\end{array}
\right.
$$
\end{example}
\begin{lemma}\label{hop}
For arbitrary small category  $\mC$ and all $n\geq 0$ 
there are isomorphisms
$\coLim_n^{\mC^{op}}\Delta_{\mC^{op}}\ZZ \cong
\coLim_n^{\mC}\Delta_{\mC}\ZZ$.
\end{lemma}
{\sc Proof.} Consider the homomorphisms
$C_n(\mC^{op},\Delta\ZZ) \rightarrow C_n(\mC,\Delta\ZZ)$, 
which assign to any element
$(c_0\stackrel{\alpha_1^{op}}\rightarrow c_1 \rightarrow \cdots 
\stackrel{\alpha_n^{op}}\rightarrow c_n)$ 
of  $C_n(\mC^{op},\Delta\ZZ)$ the element
$(c_0\stackrel{\alpha_1}\leftarrow c_1 \leftarrow \cdots 
\stackrel{\alpha_n}\leftarrow c_n)$ for even $n$, and 
 $-(c_0\stackrel{\alpha_1}\leftarrow c_1 \leftarrow \cdots 
\stackrel{\alpha_n}\leftarrow c_n)$ for odd $n$.
These homomorphisms are isomorphisms
$C_*(\mC^{op},\Delta\ZZ)\stackrel{\cong}\rightarrow C_*(\mC,\Delta\ZZ)$. 
Consequently the homology group of the complexes are 
isomorphic.
\hfill $\Box$

\medskip
\noindent
{\bf Tensor product of functors.}
Let $Lh^c: \Ab^\mC \rightarrow \Ab$ be the composition of the functors
 $L: \Set\rightarrow \Ab$ 
and $h^c(-)=\mC(c,-): \mC\rightarrow \Set$. 
For a functor $F: \mC^{op}\rightarrow \mA$ and 
an object  $A\in \mA$ we denote by  $\Hom(F,A)$
the functor $\mC \rightarrow \Ab$,
having the values  $\mA(F(c),A)$ at  $c\in Ob(\mC)$
and assigning to every  $\alpha\in \mC(c_1, c_2)$ the map 
$\Hom(F,A)(\alpha): \mA(F(c_2),A) \rightarrow \mA(F(c_1),A)$
defined as $\Hom(F,A)(\alpha)(f)=f\circ F(\alpha)$.
\begin{lemma}\label{tensor}
Let $\mC$ be a small category and $\mA$ an additive category 
with infinite coproducts. Then there exists a bifunctor 
$\otimes: \Ab^\mC\times \mA^{\mC^{op}}\rightarrow \mA$ which is additive 
at each argument and has the following properties
\begin{itemize}
\item there are 
isomorphisms
$$
      \mA(G\otimes F, A) \stackrel{\cong}\rightarrow \Ab^\mC(G, \Hom(F,A))
$$
natural in $G\in \Ab^\mC$, $F\in \mA^{\mC^{op}}$, and $A\in \mA$;
\item there are isomorphisms
$$
	Lh^c\otimes F \stackrel{\cong}\longrightarrow F(c)
$$
natural at $c\in \mC$ и $F\in \mA^{\mC^{op}}$;
\item for each  $F\in \mA^{\mC^{op}}$ the functor  
$(-)\otimes F: \Ab^{\mC} \rightarrow \mA$ preserves the colimits;
\item for every $G\in \Ab^{\mC}$ the functor
$G\otimes (-): \mA^{\mC^{op}} \rightarrow \mA$ preserves the colimits.
\end{itemize}
\end{lemma}
{\sc Proof.} Let   $\mB$ and $\mA$ be additive categories.
Denote the category of additive functors
$\mB\rightarrow \mA$ by $\Add(\mB,\mA)$.
It is well known \cite[p.11]{mit1972}, 
that for each small category $\mC$ there is the additive category 
$\Mat \ZZ\mC$ which contains $\mC$ and has the isomorphisms
 $\mA^\mC\cong \Add(\Mat\ZZ\mC, \mA)$ for every additive category $\mA$. 
There is a bifunctor \cite[Theorem 17.7.2]{sch1970} 
$$
\otimes_{\mC^{op}}: \Add(\Mat\ZZ\mC, \Ab)\times
\Add({\Mat\ZZ\mC}^{op}, \mA)
\rightarrow \mA,
$$ 
for which firstly there are natural isomorphisms
 $h^c\otimes F\cong F(c)$,
and secondly the functors $(-)\otimes F$ and $G\otimes (-)$
preserves the colimits. Using this bifunctor 
and the isomorphisms of categories
$$
\Add(\Mat\ZZ\mC, \Ab) \cong \Ab^\mC, \quad \Add({\Mat\ZZ\mC}^{op}, \mA)
\cong \mA^{\mC^{op}}
$$
we obtain the bifunctor
$\otimes: \Ab^\mC\times \mA^{\mC^{op}}\rightarrow \mA$ with 
required properties.
\hfill $\Box$ 

\medskip
\noindent
{\bf Derived functors of the colimit and partial 
derived of the tensor product at the first variable.}
Let $\mC$ be a small category.
Recall  that the homology groups of the category 
$\mC^{op}$ is denoted by  
$\coLim_n^{\mC^{op}}: \Ab^{\mC^{op}} \rightarrow \Ab$. The following 
lemma is well known \cite{obe1967} and it may be proved from general 
approach. We prove it by 
the definition  of $\coLim_n^{\mC^{op}}$ as the homology groups of the 
chain complex 
$C_*(\mC,F)$. 
\begin{lemma}\label{obe2}
For every projective resolution 
$$
	0 \leftarrow \Delta_\mC\ZZ \stackrel\varepsilon\leftarrow P_0 
 	\stackrel{d_1}\leftarrow P_1 \leftarrow \cdots 
	\stackrel{d_n}\leftarrow P_n \leftarrow \cdots
$$
the $n$-th homology group of the chain complex 
$P_*\otimes F$ is isomorphic to $\coLim_n^{\mC^{op}}F$.
\end{lemma}
{\sc Proof.} Let $P_*(\mC)$ be the chain complex in the category
 $\Ab^\mC$ where 
$$
P_n(\mC)=\bigoplus\limits_{c_0\rightarrow\cdots \rightarrow c_n} 
Lh^{c_n}\in \Ab^{\mC}
$$
are projective objects in 
 $\Ab^{\mC}$. The morphisms 
$d_n: P_n(\mC)\rightarrow P_{n-1}(\mC)$ are defined for $n>0$ as 
$\sum\limits_{i=0}^n (-1)^i d^n_i$ where 
$d^n_i: P_n(\mC)\rightarrow P_{n-1}(\mC)$ are natural transformations
which components
$(d^n_i)_c$ for $c\in \mC$
are given on generators of  
$\bigoplus\limits_{c_0\rightarrow\cdots \rightarrow c_n} 
L\mC(c_n,c)$ as 
$(d^n_i)_c(c_0\rightarrow \cdots \rightarrow c_n\rightarrow c)=
(c_0\rightarrow \cdots \rightarrow \widehat{c_i}\rightarrow \cdots
\rightarrow c_n\rightarrow c)$. Define a natural transformation 
$d_0: P_0(\mC)\rightarrow \Delta_{\mC}\ZZ$ by lettig  
$(d_0)_c(c_0\rightarrow c)=c$
on generators of the group
$P_0(\mC)(c)$ where  $\Delta_\mC\ZZ(c)$ is considered as 
the free abelian group with one generator
 $c$. 
We will prove that the sequence
$$
0\rightarrow \Delta_\mC\ZZ \stackrel{d_0}\leftarrow P_0(\mC) 
\stackrel{d_1}\leftarrow P_1(\mC) \leftarrow \cdots
\stackrel{d_n}\leftarrow P_n(\mC)\leftarrow \cdots ~.
$$
is exact in $\Ab^\mC$. 
Denote this sequence by 
 $\overline{P}_*(\mC)$. 
Let
$s_n: \overline{P}_{n-1}(\mC)(c)\rightarrow \overline{P}_{n}(\mC)(c)$ 
be homomorphisms defined for  $n>0$ on generators by
$$
s_n(c_0\rightarrow \cdots \rightarrow c_{n-1}\rightarrow c)=
(c_0\rightarrow \cdots \rightarrow c_{n-1}\rightarrow
c\stackrel{1_c}\rightarrow c),
$$
and for $n=0$ by $s_0(c)=(c\stackrel{1_c}\rightarrow c)$.
It is easy to see that $s_n$ make up the contracting homotopy of the 
chain complex $\overline{P}_*(\mC)(c)$.
Hence the sequence $\overline{P}_*(\mC)$ is exact in $\Ab^\mC$.
Thus $P_*(\mC)$ is the projective resolution of  $\Delta_\mC\ZZ$. 
It is easy to see from Lemma \ref{tensor} that the chain complex
$$
0 \leftarrow  \Delta_\mC\ZZ\otimes F \stackrel{d_0\otimes F}\leftarrow 
   P_0(\mC)\otimes F \leftarrow \cdots \leftarrow P_{n-1}(\mC)\otimes F
  \stackrel{d_n\otimes F}\leftarrow P_{n}(\mC)\otimes F \leftarrow \cdots
$$
is isomorphic to $C_*(\mC^{op},F)$. Consequently 
$H_n(P_*(\mC)\otimes F)\cong \coLim_n^{\mC^{op}}F$.

Now consider an arbitrary projective resolution 
 $P_*$ of the object 
$\Delta_\mC\ZZ$ in $\Ab^\mC$. It is homotopical equivalent to 
 $P_*(\mC)$. Since the functor $(-)\otimes F$ is additive, it 
sends homotopical equivalent chain complexes to homotopical equivalent 
chain complexes. Consequently
$H_n(P_*\otimes F)\cong H_n(P_*(\mC)\otimes F) \cong \coLim_n^{\mC^{op}}F$.
\hfill $\Box$

\medskip
\noindent
{\bf Coinitial functors.}  
A small category $\mC$ is {\em connected} if for any objects
$a,b\in \mC$ there exists a sequence 
of morphisms $a=c_0\rightarrow c_1\leftarrow c_2\rightarrow \cdots 
\leftarrow c_{2k}=b$. $\mC$ is  {\em acyclic} if
 $\coLim_n^\mC\Delta\ZZ=0$ for all $n>0$.
Let  $S: \mC \rightarrow \mD$ be a functor from 
a small category into an arbitrary category.
For any $d\in Ob(\mD)$ a {\em fibre} (or {\em a comma-category})
$S/d$ is the category which objects are pairs 
 $(c,\alpha)$, 
consisting of objects $c\in \mC$ and morphisms $\alpha\in \mD(S(c),d)$.
Morphisms in $S/d$ are defined as triples $(f, \alpha_1, \alpha_2)$ of 
$f\in \mC(c_1,c_2)$,
$\alpha_1\in \mD(S(c_1),d)$, $\alpha_2\in \mD(S(c_2),d)$
for which   $\alpha_2\circ S(f)=\alpha_1$.
A {\em forgetful functor $Q_d: S/d\rightarrow \mC$ of the fibre} 
is defined as
$Q_d(c,\alpha)=c$ on objects, and $Q_d(f, \alpha_1, \alpha_2)=f$
at morphisms.
If $S$ is a full embedding $\mC\subseteq \mD$, then $S/d$
is denoted by $\mC/d$.
\begin{definition}
A functor $S: \mC \rightarrow \mD$ between small categories is 
{\em strong coinitial} if for all object $d\in \mD$ 
the categories  $S/d$ are connected and acyclic.
\end{definition}
This is holds if and only if for all
$d\in \mD$
$$
\coLim_n^{S/d}\Delta\ZZ =\left\{
\begin{array}{cll}
0~, & \mbox{if} & n>0\\
\ZZ~, & \mbox{if} & n=0~.
\end{array}
\right.
$$
A {\em cofibre $d/S$ of an object $d\in \mD$ over 
$S: \mC\rightarrow \mD$} is a category with objects
$(c,\alpha)$ where $c\in Ob(\mC)$, $\alpha\in \mD(d, S(c))$. Its 
morphisms  
$(c_1,\alpha_1)\rightarrow (c_2,\alpha_2)$ are given by triples 
$(f,\alpha_1,\alpha_2)$ of morphisms $f\in\mC(c_1,c_2)$, 
$\alpha_1\in \mD(d, S(c_1))$, $\alpha_2\in \mD(d, S(c_2))$ satisfying 
$S(f)\circ \alpha_1 = \alpha_2$.

\begin{proposition}\label{oberst}
Let $\mC$ and $\mD$ be small categories and $S:\mC\rightarrow\mD$ a functor.
Then there are canonical homomorphisms 
$\coLim_n^{\mC^{op}}(F\circ S^{op})\rightarrow \coLim_n^{\mD^{op}}F$
which are natural in $F$.
These homomorphisms are isomorphisms for 
all $n\geq 0$ if and only if the functor  $S$ is strong coinitial.
\end{proposition}
{\sc Proof.} This follows from the assertion
 \cite[Theorem 2.3]{obe1968} affirmed for homomorphisms
$\coLim_n^{\mC}(F\circ S)\rightarrow \coLim_n^{\mD}F$. In this case 
we recive the following necessary and sufficient condition 
$\coLim_n^{d/S}\Delta\ZZ\cong \coLim_n^{pt}\Delta\ZZ$ under which 
these homomorphisms are isomorphisms.
We substitute $S$ by $S^{op}$ and obtain 
that $S$ is strong coinitial if and only if 
 $\coLim_n^{(S/d)^{op}}\Delta\ZZ\cong \coLim_n^{pt}\Delta\ZZ$. 
But  $\coLim_n^{(S/d)^{op}}\Delta\ZZ \cong \coLim_n^{S/d}\Delta\ZZ$ 
 by Lemma \ref{hop}.
\hfill $\Box$
\begin{corollary}\label{compos}
The composition of strong coinitial functors 
is strong coinitial.
\end{corollary}

\medskip
\noindent
{\bf Homology of $\mD$---sets.} 
Let   $S: \mC\rightarrow \mD$ be a functor between small categories.
Then for a cocomplete category $\mA$ the functor
$(-)\circ S: \mA^\mD\rightarrow \mA^\mC$ has (see \cite{mac1972}) 
a left adjoint functor $\Lan^S$ which value at $F\in \mA^\mC$ is called 
by a {\em left Kan extension $\Lan^S F$ of $F\in \mA^\mC$ along $S$}.
A {\em  connection components of a category} are its maximal connected 
subcategories. 
Let $\mC$ be a small category in which every connection component 
has an initial object. For each object 
$c\in \mC$ we select one initial object $i(c)$ 
in the connection component contained $c$.
Denote by $i_c: i(c)\rightarrow c$ the unique morphism 
from $i(c)$ into $c$.
Let $\init(\mC)$ be the set consisting of the selected initial objects.
\begin{lemma}\label{laninit}
Let $\mA$ be a cocomplete category.
Suppose that $S: \mC\rightarrow \mD$ is a functor between small categories
 such that for each object $d\in \mD$ every 
connection component of  $d/S$ has an initial object. 
Then for any functor
$F: \mC^{op}\rightarrow \mA$ the left Kan extension
$\Lan^{S^{op}}F: \mD^{op} \rightarrow \mA$ 
is isomorphic to a functor which values on objects
 $d\in \mD$  equal 
$\coprod\limits_{\beta \in \init(d/S)}FQ_d^{op}(\beta)$
 and assigned to every morphism $\alpha:d\rightarrow e$
of the category $\mD$ the morphism  $\overline{\alpha}$ of  $\mA$
determined by commutativity 
 $\forall \gamma\in \init(e/S)$ of the following diagrams 
$$
\begin{CD}
\coprod_{\gamma\in \init(e/S)}FQ_e^{op}(\gamma) @>\overline{\alpha}>>
\coprod_{\beta\in \init(d/S)}FQ_d^{op}(\beta) \\
@A in_\gamma AA @AA in_{i(\gamma\circ\alpha)} A  \\
FQ_e^{op}(\gamma)=FQ_d^{op}(\gamma\circ\alpha) 
@>> FQ_d^{op}(i_{\gamma\circ\alpha})  >
FQ_d^{op}(i(\gamma\circ\alpha))
\end{CD}
$$
\end{lemma}
{\sc Proof.} 
Substitute the category $\mA$ by $\mA^{op}$ in \cite[Лемма 1.1]{X1997}. 
Since products in the dual category are coproducts and 
arrows have the dual direction, the assertion follows from 
\cite[Лемма 1.1]{X1997}.
\hfill $\Box$

Let $\mD$ be a small category and 
$X\in \Set^{\mD^{op}}$ a functor. Denote by $\mD/X$ a small category 
which objects are  members
$x\in \coprod_{d\in Ob\mD}X(d)$ and its 
morphisms are triples $(\alpha, x, x')$, 
with $\alpha\in \mD(d,d')$, $x\in X(d)$, and $x'\in X(d')$ satisfying
to  $X(\alpha)(x')=x$. 

It is also concerned with $\mD$ and $X$ the fibre 
$h_*/X$ of Yoneda embedding $h_*: \mD\rightarrow \Set^{\mD^{op}}$.
Since  $h_*$ is a full embedding, every 
object $(d,\eta: h_d\rightarrow X)$
is determinent by  the second element $\eta: h_d\rightarrow X$. 
Thus  $h_*/X$ may be considered as the category with objects
$\eta: h_d\rightarrow X$ and morphisms given by commutative 
triangles
$$
  \xymatrix{
	h_d	\ar[rr] ^{\eta}
	\ar[dr]_{h_{\alpha}}
	&& X\\
	& h_{d'} \ar[ur]_{\eta'}
  }
$$
It is well known \cite{mac1972} that there is an isomorphism
 $\mD/X \stackrel{\cong}\rightarrow h_*/X$. 
It is built on objects as the inverse to the map 
$(\eta: h_d\rightarrow X)\mapsto \eta_d(1_d)$.
The value of the map at $x$  is denoted by $\widetilde{x}: h_d \rightarrow X$.

Let  $\mD$ be a small category and 
$X\in \Set^{\mD^{op}}$ and 
$F: (\mD/X)^{op} \rightarrow \Ab$ be functors.

\begin{proposition}\label{gen}
For all $n\geq 0$ there are isomorphisms 
$\coLim_n^{(\mD/X)^{op}} F \cong \coLim_n^{\mD^{op}}\breve{F}$
where  
 $\breve{F}(d)=\bigoplus\limits_{x\in X(d)}F(x)$ and 
$\breve{F}(\alpha):\breve{F}(a)\rightarrow \breve{F}(b)$ 
at $\alpha\in \mD(a,b)$ are  morphisms making the following diagrams 
commutative 
\begin{equation}\label{sq}
\begin{CD}
\bigoplus\limits_{x\in X(a)} F(x) & @>\breve{F}(\alpha)>>& 
		 \bigoplus\limits_{x\in X(b)} F(x)\\
@AA{in_x}A &  & 
		@A{in_{X(\alpha)(x)}}AA\\
F(x) & @>>F(X(\alpha)(x)\stackrel{\alpha}\rightarrow x)> & 
	F(X(\alpha)(x))\\
\end{CD}
\end{equation}
\end{proposition}
{\sc Proof.} 
Denote $\breve{F}= \Lan^{Q^{op}}F$ where $Q: h_*/X \rightarrow \mD$ is the 
forgetful full defined by the assignments  
$\widetilde{x}: h_d\rightarrow X \mapsto d\in \mD$ and  
$(\alpha, x, y) \mapsto \alpha$.
For any $a\in \mD$ each connected component of the category $a/Q$ 
has an initial object. The values of the map  
$i: Ob(a/Q)\rightarrow \init (a/Q)$  at 
$(\widetilde{x}: h_c\rightarrow X, \alpha: a\rightarrow c$ equal 
$i(\widetilde{x},\alpha)=(\widetilde{x}\circ h_\alpha, 1_a)$ and 
$i_{(\widetilde{x},\alpha)}=(\alpha,\widetilde{x}\circ
h_\alpha, \widetilde{x})$.
Since $\widetilde{x}\circ h_\alpha = \widetilde{X(\alpha)(x)}$ and 
the categories  $\mD/X$ and $h_*/X$ may be identified, 
we obtain from Lemma 
\ref{laninit}
a definition of homomorphisms  $\breve{F}(\alpha)$ as  
making squares  (\ref{sq}) commutative.

Since the coproduct functor is exact in $\Ab$, 
 $\Lan^{Q^{op}}$ is exact. Consider an arbitrary projective 
resolution of an object
$F\in \Ab^{(\mD/X)^{op}}$
$$
0 \leftarrow F \leftarrow P_0 \leftarrow P_1\leftarrow \cdots ~.
$$
The functor $\Lan^{Q^{op}}$ is left adjoint to the exact 
functor
$(-)\circ Q^{op}: \Ab^{\mD^{op}}\rightarrow \Ab^{{(\mD/X)}^{op}}$.
In accordance 
with  a dual assertion  to \cite[Prop. 6.3]{buc1972} 
it carries projective objects to projective.
Consequently it carries the projective resolution  of $F$
to the projective resolution of $\Lan^{Q^{op}}F$.
Applying the functor $\coLim^{\mD^{op}}$ to the obtained resolution we
get a sequence
$$
0\leftarrow \coLim^{(\mD/X)^{op}} \Lan^{Q^{op}} F
\leftarrow \coLim^{(\mD/X)^{op}} \Lan^{Q^{op}} P_0
\leftarrow \coLim^{(\mD/X)^{op}} \Lan^{Q^{op}} P_1 \leftarrow \cdots  
$$
such that homology groups of the chain complex  
$$
0 \leftarrow \coLim^{(\mD/X)^{op}} \Lan^{Q^{op}} P_0
\leftarrow \coLim^{(\mD/X)^{op}} \Lan^{Q^{op}} P_1 \leftarrow \cdots  
$$
are isomorphic to groups
$\coLim_n^{(\mD/X)^{op}} \Lan^{Q^{op}} F$. 
On the other hand there is the isomorphism 
$\coLim^{(\mD/X)^{op}} \Lan^{Q^{op}}\cong \coLim^{\mD^{op}}$. Thus 
the homology groups of the chain complex are isomorphic to 
 $\coLim_n^{\mD^{op}} F$.
\hfill $\Box$

\subsection{Homology of precubical sets}

\begin{definition}
A {\em  homological system of abelian groups 
on a precubical set 
$X\in \Set^{\Box_+^{op}}$} is an arbitrary functor
$F: (\Box_+/X)^{op}\rightarrow \Ab$.
\end{definition}

We will consider values of the functor
$$
\coLim^{(\Box_+/X)^{op}}_n : \Ab^{(\Box_+/X)^{op}}\rightarrow \Ab,
$$
at homological systems of abelian groups 
on precubical sets
$X: \Box_+^{op}\rightarrow \Set$.

\medskip
\noindent
{\bf Homology of a precubical abelian group.}
\begin{lemma}\label{exact}
Let $C_*$ be a chain complex
$$
0 \leftarrow L\Box_+(\II^0,\II^n) \stackrel{d_1}\leftarrow
L\Box_+(\II^1, \II^n)
\stackrel{d_2}\leftarrow \cdots 
\stackrel{d_n}\leftarrow L\Box_+(\II^n,\II^n) \leftarrow 0,
$$
consisting of free abelian groups $C_k$ generated by morphisms
$\II^k\rightarrow \II^n$. Its differentials are defined as
$d_k=\sum\limits_{i=1}^k(-1)^i (L(\partial_i^{k,1}) -
L(\partial_i^{k,0}))$ where
$\partial_i^{k,\varepsilon}(\sigma)=
\Box_+(\delta_i^{k,\varepsilon},\II^n)(\sigma)
=\sigma\circ\delta_i^{k,\varepsilon}$. Then its homology groups 
$H_q(C_*)$ are equal to $0$ if $q>0$ and $H_0(C_*)=\ZZ$.
\end{lemma}
{\sc Proof.} 
For $1\leq i\leq j\leq k-1$ 
there is equations 
$\delta_{j+1}^{k,\beta}\delta_i^{k-1,\alpha}=
 \delta_i^{k,\alpha}\delta_j^{k-1,\beta}$.
Therefore each morphism $\sigma: \II^k\rightarrow \II^n$ has an 
unique decomposition
$$
\II^k\stackrel{\delta_1^{k+1,\varepsilon_1}}\rightarrow 
\II^{k+1}\stackrel{\delta_2^{k+2,\varepsilon_2}}\rightarrow 
\II^{k+2}\rightarrow  \cdots 
\stackrel{\delta_{j_1-1}^{k+j_1-1,\varepsilon_{j_1-1}}}\longrightarrow 
\II^{k+j_1-1}
\stackrel{\delta_{j_1+1}^{k+j_1,\varepsilon_{j_1+1}}}\longrightarrow
\cdots
\stackrel{\delta_{n-k}^{n,\varepsilon_{n-k}}}\rightarrow  \II^n,
$$
which may be written as
$$
\sigma = {\delta_{n-k}^{n,\varepsilon_{n-k}}}\cdots 
{\delta_{j_1+1}^{k+j_1,\varepsilon_{j_1+1}}}
{\delta_{j_1-1}^{k+j_1-1,\varepsilon_{j_1-1}}}
\cdots {\delta_2^{k+2,\varepsilon_2}}
{\delta_{1}^{k+1,\varepsilon_{1}}}.
$$
To the studying the action of $\partial_i^{k,\varepsilon}$ at $\sigma$,
$\partial_i^{k,\varepsilon}(\sigma)=\sigma\circ\delta_i^{k,\varepsilon}$, 
we denote morphisms $\delta_i^{k,\varepsilon}$ shortly by
$\delta_i^{\varepsilon}$.
Consider a diagram:
$$
\begin{array}{lllllllllllll}
\II^k & \stackrel{\delta_1^{\varepsilon_1}}\rightarrow & \II^{k+1} &
	\stackrel{\delta_2^{\varepsilon_2}}\rightarrow  \cdots    
  \stackrel{\delta_{j_1-1}^{\varepsilon_{j_1-1}}}\rightarrow &\II^{k+j_1-1}&  
\stackrel{\delta_{j_1+1}^{\varepsilon_{j_1+1}}}\rightarrow  \cdots 
\stackrel{\delta_{j_i-1}^{\varepsilon_{j_i-1}}}\rightarrow & \II^{k+j_i}
\stackrel{\delta_{j_i+1}^{\varepsilon_{j_i}}}\rightarrow 
	\cdots
    \stackrel{\delta_{n-k}^{\varepsilon_{n-k}}}\rightarrow  \II^n\\
\uparrow\lefteqn{\delta_i^{\varepsilon}} &&  
	\uparrow\lefteqn{\delta_{i+1}^{\varepsilon}} && 
		\uparrow\lefteqn{\delta_{i+j_1-1}^{\varepsilon}} &&
		 \uparrow\lefteqn{\delta_{i+j_i-i}^{\varepsilon}}\\
\II^{k-1} &\stackrel{\delta_1^{\varepsilon_1}}\rightarrow & \II^{k} &
	\stackrel{\delta_2^{\varepsilon_2}}\rightarrow  \cdots 
	\stackrel{\delta_{j_1-1}^{\varepsilon_{j_1-1}}}\rightarrow
        &\II^{k+j_1-2}&
	\stackrel{\delta_{j_1+1}^{\varepsilon_{j_1+1}}}\rightarrow
         \cdots
	\stackrel{\delta_{j_i-1}^{\varepsilon_{j_i-1}}}\rightarrow
         & \II^{k+j_i-1}
\end{array}
$$
The lower indexes of the morphism $\delta$ increase by $1$ at the horizontal 
arrows except the arrows in the places 
 $j_1$, $j_2$, $\cdots$ where they increase by $2$. 
They increase by $1$ at the vertical arrows. 
Hence $i$ commutative squares of this diagram lead to 
the canonical decomposition of the morphism
$\sigma\circ\delta_i^{\varepsilon}$:

$$
\begin{CD}
\II^{k-1} @>\delta_1^{\varepsilon_1}>>
\II^{k} @>\delta_2^{\varepsilon_2}>> \cdots
  @>\delta_{j_1-1}^{\varepsilon_{j_1-1}}>> \II^{k+j_1-2} 
	@>\delta_{j_1+1}^{\varepsilon_{j_1+1}}>>\cdots \\
	   @. \cdots
    	    @>\delta_{j_i-1}^{\varepsilon_{j_i-1}}>> \II^{k+j_i-1} 
	      @>\delta_{j_i}^{\varepsilon}>> \II^{k+j_i} 
     @>\delta_{j_i+1}^{\varepsilon_{j_i+1}}>> \II^{k+j_i+1} @>>> \cdots
\end{CD}
$$ 
The map $\sigma$ at variables $(x_1, \cdots, x_k)$, as a composition 
of maps $\delta_1^{\varepsilon_1}$, $\delta_2^{\varepsilon_2}$, $\cdots$,
$\delta_{n-k}^{\varepsilon_{n-k}}$, has values
$$
\sigma(x_1, \cdots, x_k) = (\varepsilon_1, \cdots, \varepsilon_{j_1-1}, x_1,
 \varepsilon_{j_1+1}, \cdots, \varepsilon_{j_k-1}, x_k, \varepsilon_{j_k+1}, 
	\cdots, \varepsilon_{n-k}).
$$
Its image is equals to $I_1\times I_2 \times \cdots \times I_n$ where  
$I_{j_1}$, $I_{j_2}$, $\cdots$, $I_{j_k}$ are nondegenerated intervals.
The image of the map $\sigma\circ\delta_i^{\varepsilon}$ is equal to 
$$
I_1\times \cdots \times I_{j_i-1}\times \{\varepsilon\}\times 
I_{j_i+1}\times \cdots \times I_n.
$$
The image $I_1\times \cdots \times I_n$ defines the map $\sigma$.
Thus the differentials of the precubical set 
correspond to the differential 
of the chain complex
$C_*(\mK([0,1]^m))$.
Using Lemma \ref{KMM} we obtain the required assertion.
\hfill $\Box$
\begin{remark}\label{contra}
Similar assertion remains true if we substitute 
$\Box_+$  by a category $\Delta$ of finite totally ordered 
sets 
$[n]=\{0, 1, \cdots, n\}$, $n\geq 0$ and nondecreasing maps.
Differentials  
$d_k: L\Delta([k],[n])\rightarrow L\Delta([k-1],[n])$ are defined as
$d_k(\sigma)=\sum\limits_{i=0}^k (-1)^i L\Delta(\delta^k_i,[n])(\sigma)$
where $\delta^k_i: [k-1]\rightarrow [k]$ is increasing 
maps such that $i\not\in \Imm(\delta^k_i)$.
Similar assertion also is true for a category 
 $\Delta_+$ of finite totally ordered sets $[n]$, $n\geq 0$,
and increasing maps \cite[Lemma 1.2]{X1997}. But it is false for  
a category $\Box$ with  
objects $Ob(\Box)=Ob(\Box_+)$ and morphisms generated 
by  $\delta^{n,\varepsilon}_i$ 
and nondecreasing surjections 
$\sigma^n_i: \II^n\rightarrow \II^{n-1}$, 
$\sigma^n_i(x_1, \cdots, x_n)=
(x_1, \cdots, x_{i-1}, x_{i+1}, \cdots, x_n )$.
Take for example  $n=0$. Then we obtain a chain complex
$$
0\leftarrow \ZZ \stackrel{0}\longleftarrow 
\ZZ \stackrel{0}\longleftarrow \ZZ \stackrel{0}\longleftarrow \cdots
$$
which all homology groups are isomorphic to
 $\ZZ$. 
\end{remark}
A prototype of the following proposition is 
 \cite[Lemma 4.2]{gab1967} which is proved there for simplicial objects 
of an abelian catagory.  We prove it for precubical abelian groups.
\begin{proposition}\label{hlim}
Let $F:\Box_+^{op} \rightarrow \Ab$ be a precubical abelian group
and 
$F_*$ a chain complex  $F_n= F(\II^n)$ with differentials 
$d_n = \sum\limits_{i=1}^n (-1)^i(F(\delta^{n,1}_i)-F(\delta^{n,0}_i))$.
Then for all $n\geq 0$ 
there are  isomorphisms
$\coLim_n^{\Box_+^{op}} F \cong H_n(F_*)$
natural in $F\in \Ab^{\Box_+^{op}}$.
\end{proposition}
{\sc Proof.} It follows from Lemma \ref{exact} that the sequence 
$$
0\leftarrow \Delta_{{\Box}_+}\ZZ \stackrel{d_0}\leftarrow Lh^{\II^0} 
	\stackrel{d_1}\leftarrow Lh^{\II^1} \leftarrow \dots 
		\stackrel{d_k}\leftarrow Lh^{\II^k} 
	\stackrel{d_{k+1}}\leftarrow \dots
$$
is the projective resolution of $\Delta_{{\Box}_+}\ZZ\in Ab^{{\Box}_+}$.
We have  isomorphisms  by Lemma \ref{obe2}
$\coLim_n F \cong  H_n(P_*\otimes F)$.
It follows from Lemma \ref{tensor} that
there are natural isomorphisms
$Lh^{\II^n}\otimes F \stackrel{\cong}\rightarrow F(\II^n)$.
The differentials of the chain complex 
$Lh^{\II^*}\otimes F$ are equal to   
$\sum\limits_{i=1}^n (-1)^i(Lh^{\delta^{n,1}_i}-
Lh^{\delta^{n,0}_i})\otimes F$.
The isomorphism  
$Lh^{(-)}\otimes F \stackrel{\cong}\rightarrow F(-)$  
carries these differentials to homomorphisms  
$\sum\limits_{i=1}^n (-1)^i(F(\delta^{n,1}_i)-F(\delta^{n,0}_i))$. 
From this we obtain isomorphisms 
$\coLim_n^{{\Box}_+^{op}} F \cong H_n(F_*)$ which are natural in $F$.
\hfill $\Box$

In accordance with Remark \ref{contra} this proposition is true 
in the cases of the categories $\Delta$ and $\Delta_+$.
The category $\Box$ has the terminal object $\II^0$. 
Hence the functor $\coLim^{\Box}$ is exact and we obtain by Lemma
 \ref{hop} that   
$\coLim_n^{\Box^{op}}\Delta\ZZ\cong \coLim_n^{\Box}\Delta\ZZ = 0$ for $n>0$.
We saw that  $H_n(\Delta\ZZ_*)\cong \ZZ$ for all $n\in \NN$ in this case.
Consequently this proposition is false for the category $\Box$.

\medskip
\noindent
{\bf Homology groups with coefficients in a homological system.} 
Now we study homology groups
 $\coLim_n^{(\Box_+/X)^{op}}F$ 
of a precubical set  $X\in \Set^{\Box_+^{op}}$ with coefficients 
in a homological system $F$ on $X$. 
Consider abelian groups
$C_n(X,F)=\bigoplus\limits_{\sigma\in X_n}F(\sigma)$. Define differentials
$d_i^{n,\varepsilon}: C_n(X,F)\rightarrow C_{n-1}(X,F)$ as 
homomorphisms making the following diagrams commutative 
$$
\begin{CD}
\bigoplus\limits_{\sigma\in X_n}F(\sigma) @>d_i^{n,\varepsilon}>> 
		\bigoplus\limits_{\sigma\in X_{n-1}}F(\sigma)\\
@A in_\sigma AA @A in_{\sigma\delta^{\varepsilon}_i}AA\\
F(\sigma) 
@>F(\delta_i^{\varepsilon}:
\sigma\delta^{\varepsilon}_i\rightarrow \sigma)>>
F(\sigma\delta_i^{\varepsilon})
\end{CD}
$$

\begin{definition}
Let   $F: (\Box_+/X)^{op}\rightarrow \Ab$
be a homological system of abelian groups on 
a precubical set  $X$. 
{\em Homology groups   $H_n(X,F)$ with coefficients in $F$}
are $n$-th homology groups of the chain complex $C_*(X,F)$ consisting 
of abelian groups  $C_n(X,F)= \bigoplus\limits_{\sigma\in X_n}F(\sigma)$ 
and differentials $d_n=\sum\limits_{i=1}^n (-1)^i (d^{n,1}_i-d^{n,0}_i)$. 
\end{definition}

\begin{example}\label{gou3}
Consider functors $\ZZ^0, \ZZ^1: {\Box}_+^{op}\rightarrow \Ab$ 
with values  $\ZZ^{\alpha}(\II^n)=\ZZ$ at objects for all $n\in\NN$.
Its values at morphisms are defined by homomorphisms
$\ZZ^{\alpha}(\delta_i^{n,\varepsilon}): \ZZ\rightarrow \ZZ$ 
which equal $1_{\ZZ}$ for $\varepsilon=\alpha$ and  $0$ for
$\varepsilon\not=\alpha$. For any precubical set
$X= (X_n, \partial_i^{n,\varepsilon})$ are defined 
homological systems
$\ZZ^0\circ Q_X^{op}, \ZZ^1\circ Q_X^{op}:
 ({\Box}_+/X)^{op}\rightarrow \Ab$
where  $Q_X: {\Box}_+/X \rightarrow {\Box}_+$ is the forgetful 
functor of a fibre.
We get two chain complexes $C_*(X,\ZZ^0 \circ Q_X^{op})$ and 
$C_*(X,\ZZ^1 \circ Q_X^{op})$. 
Complexes $C_*(X,\ZZ^\alpha\circ Q_X^{op})$ for $\alpha\in \{0,1\}$ 
consist of abelian groups
$$
0\leftarrow L(X_0) \stackrel{d_1^\alpha}\leftarrow L(X_1)
\stackrel{d_2^\alpha}\leftarrow L(X_2) \leftarrow \dots
$$
with differential
$d_n^0= -\sum\limits_{i=1}^n (-1)^i L(\partial_i^{n,\varepsilon})$ 
and 
$d_n^1= \sum\limits_{i=1}^n (-1)^i L(\partial_i^{n,\varepsilon})$. 
Therefore every of two types 
homology groups  
considered by Goubault \cite{gou1995}
are particular cases of the homology 
groups with coefficients in a homological system.
\end{example}

Let  $\mHS$ be a category each object of which is a pair 
$(X,F)$ consisting of 
a precubical set $X$ and homological system 
$F: (\Box_+/X)^{op} \rightarrow \Ab$. Morphisms $(X,F)\rightarrow (Y,G)$
are given by morphisms of precubical sets
 $f: X\rightarrow Y$ for which the following diagram is commutative 
$$
  \xymatrix{
	(\Box_+/X)^{op}	\ar[rr] ^{F}
	\ar[dr]_{(\Box_+/f)^{op}}
	&& \Ab\\
	& (\Box_+/Y)^{op} \ar[ur]_{G}
  }
$$
Here $(\Box_+/f)^{op}$ denotes a functor dual
to $(\Box_+/f)(\sigma)= f\circ \sigma$. 
Since the correspondence $(X,F)\mapsto \breve{F}$ give  
the functor $\breve{(-)}: \mHS \rightarrow \Ab^{\Box_+^{op}}$, we obtain
functors $(X,F)\mapsto H_n(X,F)$
and $(X,F)\mapsto \coLim_n^{(\Box_+/X)^{op}}{F}$.

\begin{theorem}\label{comcub}
There are isomorphisms
$H_n(X,F) \cong \coLim_n^{(\Box_+/X)^{op}}F$ 
which are natural in $(X,F)\in \mHS$
 for all $n\geq 0$.
\end{theorem}
{\sc Proof.} 
We have 
 $\coLim_n^{(\Box_+/X)^{op}}F \cong $ $ \coLim_n^{\Box_+^{op}} \breve{F}$ 
by Proposition \ref{gen}.
It follows from the property of Kan extensions  
$\Lan^S\circ \Lan^T \cong \Lan^{S\circ T}$ that this isomorphism 
is natural in  $(X,F)\in \mHS$.
The chain complex  $\breve{F}_*$ is natural isomorphic 
to $C_*(X,F)$.
It follows from Proposition \ref{hlim} that 
$\coLim_n^{\Box_+^{op}} \breve{F}\cong H_n(\breve{F}_*)$.
The isomorphisms
$\coLim_n^{\Box_+^{op}} G \stackrel{\cong}\rightarrow H_n(G_*)$ 
 are natural in precubical abelian groups
$G\in \Ab^{\Box_+^{op}}$. Therefore we have the natural isomorphisms
$\coLim_n^{\Box_+^{op}} \breve{F}
\stackrel{\cong}\rightarrow H_n(\breve{F}_*)$.
Consequently there is a natural isomorphism 
 $H_n(C_*(X,F))\cong \coLim_n^{(\Box_+/X)^{op}}F$.
\hfill $\Box$

\begin{example}\label{concub}
Let ${\Box}^0=h_{\II^0}\in \Set^{{\Box}^{op}}$ be a cubical point. 
The category  $\Box/\Box^0$ is isomorphic to $\Box$. Hence
  $\coLim_n^{({\Box}/{\Box}^0)^{op}}\Delta\ZZ\cong \coLim_n^{pt}\Delta\ZZ$.
But $H_n({\Box}^0,\Delta\ZZ)\cong \ZZ$ for all $n\geq 0$.
Consequently Theorem \ref{comcub} is false for 
cubical sets.
\end{example}

\medskip
\noindent
{\bf A criterion of a homological isomorphism.} 
Let $f: X\rightarrow Y$ be a morphism of precubical sets.
We identify cubes $\sigma\in X_m$ with corresponding natural 
transformations
 $\widetilde{\sigma}: h_{\II^m} \rightarrow X$ which 
 are called by {\em singular cubes}.
The functor $\Box_+/f: \Box_+/X \rightarrow \Box_+/Y$ carries 
every cube $\sigma: h_{\II^m}\rightarrow X$
to the composition $f\circ\sigma$.
It assigns to any morphism $(\delta,\sigma,\tau)$ of the category $\Box_+/X$   
the morphism $(\delta,f\circ\sigma,f\circ\tau)$.
Let  $f^*: \Ab^{(\Box_+/Y)^{op}} \rightarrow \Ab^{(\Box_+/X)^{op}}$ 
be a functor assigning to each homological system 
$F$ on $Y$ the composition 
  $F\circ (\Box_+/f)^{op}$. The functor $f^*$ acts at morphisms 
$\eta: F_1\rightarrow F_2$ by Godement's law
$f^*(\eta)=\eta*((\Box_+/f)^{op})$, so 
$f^*(\eta)_\sigma= \eta_{f\circ\sigma}$.
It follows from Proposition \ref{oberst} that there are canonical 
  homomorphisms 
$\coLim_n^{(\Box_+/X)^{op}}f^*(F) \rightarrow \coLim_n^{(\Box_+/Y)^{op}}F$
which are natural in $F$.
Thus we have the homomorphisms
$H_n(X, f^*(F))\rightarrow H_n(Y,F)$. 

Let $\overleftarrow{f}(\sigma)$ be a precubical set defined by 
a pullback diagram
\begin{equation}\label{invimage}
\begin{CD}
X @>f>> Y\\
@A f_{\sigma}AA   @AA \sigma A\\
\overleftarrow{f}(\sigma) @>>> h_{\II^n}\\
\end{CD}
\end{equation}
The precubical set  $\overleftarrow{f}(\sigma)$ is 
an {\em inverse image of the singular cube  $\sigma$}.
Since $(\Box_+/f)/\sigma \cong \Box_+/\overleftarrow{f}(\sigma)$ and  
$\coLim_n^{(\Box_+/\overleftarrow{f}(\sigma))^{op}}\Delta\ZZ \cong 
H_n(\overleftarrow{f}(\sigma), \Delta\ZZ)$, 
it follows from Lemma \ref{hop} and Proposition  \ref{oberst} the following
\begin{corollary}
The homomorphisms $H_n(X, f^*(F))\rightarrow H_n(Y,F)$ 
are isomorphisms at all 
$F\in \Ab^{(\Box_+/Y)^{op}}$ and $n\geq 0$ if and only if 
  $H_n(\overleftarrow{f}(\sigma),\Delta\ZZ)$ are isomorphic to 
$\coLim_n^{pt}\Delta\ZZ$
for all   $\sigma\in \Box_+/Y$ and $n\geq 0$. 
\end{corollary}

\medskip
\noindent
{\bf Homology groups of a colimit of precubical sets.} 
Let  $\mA$ be an arbitrary category. Functors from small categories into 
 $\mA$ are called  
{\em diagrams of object of the category  $\mA$}.
If a domain of a diagram is a small category
 $J$, then this diagram is denoted by $\{X^i\}_{i\in J}$ where $X^i$ 
are its values at $i\in Ob(J)$. In particular, 
{\em a diagram of precubical sets }
$\{X^i\}_{i\in J}$ consists of precubical sets $X^i$ and morphisms 
 $X^\alpha: X^i\rightarrow X^j$ which are given for  
 $\alpha\in I(i,j)$ and satisfy to equations
 $X^{\alpha\circ\beta}= X^{\alpha}\circ X^{\beta}$ and  
$X^{1_i}=1_{X^i}$.

Recall  that a first quadrant spectral sequence 
 $E^r_{p,q}$ \cite{mac1963} is said to
 {\em converge} to a graded abelian group  
$\{H_n\}$, and we write in this case 
$E^r_{p,q} {~~}_p \Rightarrow H_{p+q}~$
if for each  $n\in \NN$ there is a sequence 
of subgroups
$$
0=F_{-1}H_n \subseteq F_{0}H_n
\subseteq F_{1}H_n \subseteq \cdots \subseteq F_n H_n = H_n,
$$
such that 
$F_p H_n / F_{p-1}H_n \cong E^{\infty}_{p,n-p}$ for all $0\leq p\leq n$.
\begin{proposition}
Let  $\{X^i\}_{i\in J}$ be a diagram of precubical sets 
$X^i\in \Set^{\Box_+^{op}}$ such that
\begin{equation}\label{cond}
\coLim_q^J\{L(X^i_n)\}_{i\in J}= 0, \mbox{ for any } n\in \NN \mbox{ and }
q>0.
\end{equation}
Let $\lambda_i: X^i\rightarrow
\coLim^J\{X^i\}_{i\in J}$ is the colimit cone.
Then for every  homological system
$$
F: (\Box_+/\coLim^J\{X^i\}_{i\in J})^{op} \rightarrow \Ab
$$ 
there exists a first quadrant spectral sequence
with  
$E^2_{p,q}= \coLim_p^{J}\{H_q(X^i,\lambda_i^*(F))\}_{i\in J}$
which converges to graded abelian group
$\{H_n(\coLim^J\{X^i\}_{i\in J}, F)\}$.
\end{proposition}
{\sc Proof.} Let $\mA$ be an abelian category with exact products 
and $F: \Box_+/\coLim^J\{X^i\}_{i\in J} \rightarrow \mA$ a functor.
Under the condition (\ref{cond}) by 
\cite[Cor. 2.4]{X1989} there is a spectral sequence 
$$
E_2^{p,q}= \lim\nolimits^p_{J^{op}} 
\{ \lim\nolimits^q_{\Box_+/X^i} F\circ (\Box_+/\lambda_i) \}_{i\in J}
 {~~}_p\Rightarrow \lim\nolimits^{p+q}_{\Box_+/\coLim^J\{X^i\}_{i\in J}} F
$$
Substitute  $\mA = \Ab^{op}$. We have obtained 
a spectral sequence 
$$
E^2_{p,q}=
\coLim_p^{J}\{\coLim_q^{(\Box_+/X^i)^{op}}\lambda_i^*(F)\}_{i\in J}
{~~}_p\Rightarrow \{\coLim_{p+q}^{(\Box_+/\coLim^J\{X^i\})^{op}}F\}.
$$
Theorem \ref{comcub} leads us to the wanted spectral sequence.
\hfill $\Box$

\begin{definition}
Let  $X$ be a precubical set. {\em A locally directed cover 
of  $X$}  is a diagram of precubical sets $\{X^i\}_{i\in J}$ 
over a partially ordered set 
$J$ satisfying the following three conditions

1) it assigns to each pair of members $i\leq j$ of the set $J$
an inclusion $X^i \subseteq X^j$;

2) for each $n\in \NN$ the equation $X_n= \bigcup_{i\in J}X_n^i$ holds true;

3) for each $\sigma\in X_n^i\cap X_n^j$ there exists $k\in J$ such that 
$k\leq i$ and $k\leq j$ and $\sigma\in X_n^k$.
\end{definition}

\begin{corollary}\label{locdir}
Let $\{X^i\}_{i\in J}$ be a locally directed
cover of a precubical set $X$ with
inclusions $\lambda_i$ inclusions $X^i \subseteq X$ and 
$F$ a homological system of abelian groups 
on $X$.
Then there exists a first quadrant spectral sequence with 
$E^2_{p,q}= \coLim_p^{J}\{H_q(X^i,\lambda_i^*(F))\}_{i\in J}$ 
which  converges  
to $\{H_n(X, F)\}$.
\end{corollary}
{\sc Proof.} It follows from Theorem \ref{comcub} and 
\cite[Cor. 2.5]{X1989}.
\hfill $\Box$
\begin{example}
Consider a cover of a cubical subset $X$ of Euclidian space by the 
elementary cubes. These cubes put together a partially ordered aet 
 $J$ ordered by the inclusion.
The diagram of these cubes and inclusions is 
the locally directed cover. 
For the homological system $F= \Delta\ZZ$ 
the spectral sequence collapses and gives  
 isomorphisms
$H_n(X,\Delta\ZZ)\cong \coLim_n^J\Delta\ZZ$ для всех $n\in \NN$.
\end{example}

\medskip
\noindent
{\bf A spectral sequence of a morphism.} Let $f: X\rightarrow Y$ be a
morphism of cubical sets. For each $\sigma\in \Box_+/Y$ the pullback
diagram (\ref{invimage}) defines the inverse image of the cube
 $\sigma$ with the morphism  
$f_{\sigma}: \overleftarrow{f}(\sigma) \rightarrow X$.
The inverse images of the singular cubes make up a diagram 
$\{\overleftarrow{f}(\sigma)\}_{\sigma\in \Box_+/Y}$ with the colimit 
isomorphous to $X$. Using a spectral sequence of a morphism  
\cite[Теорема 4.1]{X1991} where we substitute   $\mA$ by the 
category  $\Ab^{op}$ obtain the following assertion by 
Theorem \ref{comcub}:
\begin{corollary}\label{spmor}
Let $f: X\rightarrow Y$ be a morphism of precubical sets and 
$F: (\Box_+/X)^{op}\rightarrow \Ab$ homological system 
of abelian groups on $X$. The there exists 
a first quadrant spectral sequence 
$$
E^2_{p,q} = \coLim_p^{\Box_+/Y}\{H_q(\overleftarrow{f}(\sigma),
f_{\sigma}^*(F))\}_{\sigma\in{\Box_+/Y}} {~~}_p\Rightarrow 
H_{p+q}(X,F).
$$
\end{corollary}
In particular, if the diagram 
$\{H_q(\overleftarrow{f}(\sigma) 
f_{\sigma}^*(F))\}_{\sigma\in{\Box_+/Y}}$ 
consists of isomorphisms, then we can to obtain 
a diagram over the dual category by inverting the isomorphisms. Denote it by
 $\{H_q(\overleftarrow{f}(\sigma),
f_{\sigma}^*(F))\}_{\sigma\in(\Box_+/Y)^{op}}^{-1}$. Using
Theorem \ref{comcub} we get the spectral sequence
$$
E^2_{p,q}\cong H_p(Y, \{H_q(\overleftarrow{f}(\sigma),
f_{\sigma}^*(F))\}_{\sigma\in(\Box_+/Y)^{op}}^{-1}) {~~}_p\Rightarrow 
H_{p+q}(X,F)~,
$$
concerning homology of  $X$ and  $Y$.

\subsection{The category of factorizations}

Let $\mC$ be a small category. Denote by $\fF\mC$ the {\em category of factorization in $\mC$} 
\cite{bau1985}. Recall that $\Ob(\fF\mC)= {\Mor}(\mC)$ and for every $\alpha, \beta \in \Ob(\fF\mC)$
the set of morphisms $\fF\mC(\alpha,\beta)$ consists of the pairs
 $(f,g)$  of morphisms $f, g \in \Mor\mC$ 
$$
  \xymatrix{
	\circ	\ar[rr] ^{g}
	&& \circ \\
	\circ \ar[u]_{\alpha}
	&&  \ar[ll]_{f} \circ \ar[u]_{\beta}
  }
$$
satisfying
$g\circ\alpha\circ f = \beta$. The composition of
 $\alpha\stackrel{(f_1,g_1)}\rightarrow \beta$ and 
$\beta\stackrel{(f_2,g_2)}\rightarrow \gamma$
 is defined by 
$\alpha\stackrel{(f_1\circ f_2,g_2\circ g_1)}\longrightarrow \gamma$.
The identity morphism of an object $a\stackrel\alpha\rightarrow b$ in  
$\fF\mC$ is the pair of the identity morphisms
$ \alpha \stackrel{(1_a, 1_b)}\longrightarrow \alpha$.

\begin{lemma}\label{factop}
For any category $\mC$ there is an isomorphism $\fF(\mC^{op})\cong \fF\mC$.
\end{lemma}

\begin{corollary}\label{factopcor}
Let $\mC$ be a category and $\alpha$ its morphism. 
 $\fF\mC/\alpha$ is isomorphic to the category with objects triples
of morphisms $(x, \beta, y)$ for which the composition is defined and 
$x\circ\beta\circ y = \alpha$ and with morphisms $(x, \beta, y)\rightarrow (x', \beta', y')$ 
commutative diagrams
$$
	\xymatrix{
		& \circ \ar[ld]_x \ar[dd]_v & \circ \ar[l]^\beta 	\\
		\circ &&& \circ \ar[lu]_y \ar[ld]^{y'}\\
		& \circ \ar[lu]^{x'} & \circ \ar[l]_{\beta'} \ar[uu]_w
	}
$$
\end{corollary}

In \cite{bau1985} it was introduced the cohomology of categories with 
coefficients in natural system. It follows from \cite[Theorem 4.4]{bau1985} that 
this cohomology may be defined as the right derived functor $\Lim^n_{\fF\mC}$ of limit.
We give a dual definition.

\begin{definition}
Let $\mC$ be a small category and $\fF \mC$ the category of factorizations. 
A {\em contravariant natural system on $\mC$} is a functor  
$F: (\fF \mC)^{op}\rightarrow \Ab$. 
For $n\geq 0$ the {\em $n$-th homology group of $\mC$ with coefficients in
$F$} is the abelian groups $\coLim_n^{(\fF\mC)^{op}}F$.
\end{definition}

Since the category of factorization has a complicated structure, we seek its strong 
coinitial subcategories with simple building.  To this end we will study comma 
categories $\fF\mC/\alpha$. 

A category is said to be {\em cancellative} if its morphisms are epimorphic and monomorphic.
A category $\mC$ is cancellative if and only if the implications 
$\alpha\gamma=\beta\gamma~ \Rightarrow~ \alpha=\beta$ and 
 $\gamma\alpha=\gamma\beta~ \Rightarrow~ \alpha=\beta$ are true for all $\alpha, \beta, \gamma \in \mC$.
The following lemma is useful for cancellative categories.

Recall that a preordered set may be defined as a small category
$\mC$ such that  for any $a,b\in \Ob\mC$ the set  $\mC(a,b)$ is either empty or has 
precisely one element. In the letter case we write $a\leq b$.
If moreover $a\leq b$ and $b\leq a$ implies $a=b$, then $\mC$ is the partially ordered set.

\begin{lemma}\label{commapart}
Let $\mC$ be a small cancellative category. If $\mC$ does not contain nonidentity retractions,
 then  the comma category $\fF\mC/\alpha$ is the partially 
ordered set for each $\alpha\in Ob(\fF\mC)$.
\end{lemma}
{\sc Proof.} For any $\alpha \in \Mor\mC$ the objects of the category $\fF\mC/\alpha$ 
are triples of morphisms $(x,w,y)$  of $\mC$ for which $y\circ w\circ x= \alpha$.
 Consider any (parallel) pair of morphisms
$(x,w,y){{(u_1,v_1)\atop\longrightarrow}\atop{\longrightarrow\atop(u_2,v_2)}}(x',w',y')$
of the category $\fF\mC/\alpha$.
These morphisms are given by the commutative diagram 
$$
  
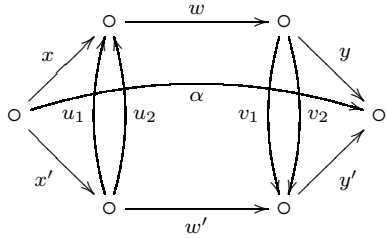
\xymatrix{
	& \circ	\ar[rr] ^{w}  & & \circ \ar[dr] ^{y} \ar@/^/[dd]^{v_2} \ar@/_/[dd]_{v_1}\\
	\circ \ar[ur]^x \ar@/^1pc/[rrrr] _\alpha  \ar[dr]_{x'} &&&& \circ \\
	& \circ \ar[rr]_{w'} \ar@/^/[uu]^{u_1} \ar@/_/[uu]_{u_2} && \circ \ar[ur]_{y'}
  }
$$
It follows from $y'v_1= y= y'v_2$ that $v_1= v_2$. Similarly,  $u_1 x'= x= u_2 x'$ implies $u_1= u_2$. 
Hence $u_1= u_2$ and $v_1= v_2$. Consequently, any two morphisms 
$(x,w,y)\rightarrow (x',w',y')$ are same. It follows that $\fF\mC/\alpha$ is the preordered set.
If $\mC$ does not contain nonidentity retraction, then any two morphisms
$$
(x,w,y){{(u_1,v_1)\atop\longrightarrow}\atop{\longleftarrow\atop(u_2,v_2)}}(x',w',y')
$$
lead to the equalities  $u_1 x'= x$, $u_2 x=x'$, $y'v_1= y$, $y v_2 = y'$. It follows that 
 $u_1 u_2 x = x$, $y v_2 v_1= y$. Since $\mC$ is cancellative, we obtain that 
 $u_1 u_2$ and $v_2 v_1$ are identity morphisms. If $\mC$ does not contain 
 nonidentity retraction, then it results from this that $u_1$, $u_2$, $v_1$, $v_2$ are 
identity morphisms. Therefore the relation $\leq$ is antisymmetric and 
 $\fF\mC/\alpha$ is the partially ordered set.
\hfill $\Box$

\subsection{Acyclic partially ordered sets}

Recall that a small category $\mC$ is {\em acyclic} 
if $H_n(\mC)\cong H_n(pt)$ for all $n\geq 0$.
The study of strong coinitial functors is closely related with the acyclic property.
We point out two simple conditions of the acyclic property of a partially ordered set.

\begin{lemma}\label{covlem}
Let $X$ be a partially ordered set. Let $V$ and  $W$ be its closed subsets for which
 $X=V\cup W$. If
$V$, $W$, and  $V\cap W$ are acyclic, 
then $X$ is acyclic.
\end{lemma}
{\sc Proof.} Let $C_*(X)=C_*(X,\Delta_X\ZZ)$. The unique functor $X\rightarrow pt$ gives the 
chain homomorphism 
$C_*(X)\rightarrow C_*(pt)$. Denote by $\tilde{C}_*(X)$ its kernel.
It is easy to see that $X$ is acyclic if and only if 
 $H_n(\tilde{C}_*(X))=0$ for all $q\geq 0$. 
The subsets $V^{op}$ and $W^{op}$ are open and
 $X^{op}= V^{op}\cup W^{op}$. 
There exists a exact sequence 
$$
0\rightarrow \tilde{C}_*(V^{op}\cap W^{op})\rightarrow \tilde{C}_*(V^{op})\oplus\tilde{C}_*(W^{op})
\rightarrow \tilde{C}_*(X^{op})\rightarrow 0.
$$ 
A corresponding long exact sequence 
leads us to isomorphisms
 $H_n(\tilde{C}_*(X^{op}))=0$ for all $q\geq 0$.
\hfill $\Box$

An easy induction proves

\begin{corollary}\label{covcor}
Let $X = \bigcup\limits_{i=1}^n W_i$ be the union of closed subsets $W_i\subseteq X$.
If  the intersections
 $W_{i_1}\cap W_{i_2}\cap \cdots \cap W_{i_k}$  are acyclic for all 
$\{i_1, i_2, \cdots ,i_k\} \subseteq \{1, 2, \cdots , n\}$, 
then $X$ is acyclic.
\end{corollary}

\section{Homology groups of free partially commutative monoids}

From here  $n$-th integral homology groups  $\coLim_n^{\mC}\Delta\ZZ$ 
of a small category $\mC$ will be denoted by $H_n(\mC)$.
At first we consider the case of commutative monoids 
and then generalize obtained assertions to partially commutative monoids.

\subsection{Homology of a commutative monoid}

{\bf Factorization category of the free commutative monoid.} 
\begin{definition}
A {\em $n$-th Leech homology group   $H_n(M,F)$ with coefficients 
in a contravariant natural system
$F: (\fF M)^{op}\rightarrow \Ab$} is a group $\coLim_n^{(\fF M)^{op}}F$.
\end{definition}

\begin{lemma}\label{strcof}
Let  $\NN=\{1, a, a^2, \ldots\}$ be a free monoid with one generator
$a$ and $T=\{1,a\}$ its subset and  $\fF{T}\subset \fF\NN$ a full 
subcategory of the factorization category of $M$ with the set of objects 
$T$. Then the inclusion $\fF{T}\subset \fF\NN$ is strong coinitial.
\end{lemma}
{\sc Proof.} 
The category  $\fF{T}$ consists of two objects $1$ and $a$
and two morphisms 
$1\stackrel{(1,a)}\rightarrow a$ and $1\stackrel{(a,1)}\rightarrow a$ 
except the identities.  
For any  $p\in \NN$ objects of  
$\fF{T}/a^p$ may be partioned of the following two classes
\begin{itemize}
\item $1\stackrel{(a^s,a^t)}\longrightarrow a^p$ where $s+t=p$, 
for $p\geq 0$,
\item $a\stackrel{(a^s,a^t)}\longrightarrow a^p$ where $s+t=p-1$, 
for $p\geq 1$. 
\end{itemize}
The fibre $\fF{T}/a^p$ does not contain 
nonidentity morphisms between objects of the same classes.
Between any two objects there exists one morphism at the most.
Every nonidentity morphism of the fibre acts from an object of 
the first class to an object of the second class and decreases a 
degree by 
 $1$ either a first or second morphism. Therefore the category
$\fF{T}/a^p$ is isomorphic to a category 
consisting of the morphisms
$$
(1,a^p)\rightarrow (1,a^{p-1})\leftarrow (a^1, a^{p-1}) \rightarrow \cdots
$$
$$
~~~~~~~~~~~~~~~~~~~
\cdots \leftarrow (a^s, a^t) \rightarrow (a^s, a^{t-1})
\leftarrow (a^{s+1}, a^{t-1})
\rightarrow \cdots \leftarrow (a^p,1)
$$

It is well known (for example, see \cite[Prop. 2.2]{X1997} for the dual 
assertion) that homology groups of any partially ordered set
 $(P,\leq)$ with coefficients in a functor $G: P\rightarrow \Ab$ 
are isomorphic to homology groups of a subcomplex $C_*^+(P,G)$ of
$C_*(P,G)$
(see Definition \ref{defhomolcat}) which consists of 
abelian groups $C_n^+(P,G)=\bigoplus_{p_0<\cdots <p_n} G(p_0)$ where 
$p_0, \cdots, p_n$ 
are members
of $P$ for which the strong inequalities $p_0< \cdots < p_n$ hold.
If $G=\Delta\ZZ$ then the subcomplex consists of 
free abelian groups
 $\bigoplus_{p_0<\cdots <p_n} \ZZ$ and differentials 
defined on generators by
$d_n ({p_0<\cdots <p_n}) = 
\sum\limits_{i=0}^n (-1)^i({p_0<\cdots<\widehat{p_i} <\cdots <p_n})$.
From here one can obtain the following chain complex 
for the calculating the integral homology groups 
of the partial set $\fF{T}/a^p$ 
$$
0\rightarrow \ZZ^{2p}\stackrel{d_1}\longrightarrow \ZZ^{2p+1}\rightarrow 0,
$$
The differential $d_1$ has the matrix
$$
D_1 = 
\left(
\begin{array}{rrrrrrr}
-1  & 0 &  0 &\cdots & 0 & 0\\
1 & 1& 0 &\cdots & 0 & 0\\
0 & -1& -1 &\cdots & 0 & 0\\
\vdots& \vdots& \vdots & \ddots & \vdots& \vdots\\
0  & 0 & 0&\cdots & 1 & 1\\
0  & 0 & 0&\cdots & 0 & -1
\end{array}
\right)
$$
consisting of  $2p+1$ strings and $2p$ columns.
Since for every column $x\in \ZZ^{2p}$ the equation  $D_1 x=0$
has the unique solution $x=0$, the homomorphism $d_1$ is injective. 
Hence the homology groups equal 
$H_1({\fF{T}/a^p})=0$, 
$H_0({\fF{T}/a^p})=\ZZ$.
\hfill $\Box$

Consider the free commutative monoid
$\NN^n$ generated 
by $a_1$, $a_2$, $\cdots$, $a_n$.
Let $T^n\subset \NN^n$ be a subset consisting of finite 
productions  $a_{i_1}a_{i_2}\cdots a_{i_k}$ for which 
$1\leq i_1< i_2 < \cdots < i_k\leq n$. 
The category $\fF\NN^n$ contains a full subcategory with a set of 
objects $T^n$. Denote this subcategory  by 
 $\fF T^n$.

\begin{lemma}\label{confn}
The inclusion $\fF T^n \subset \fF\NN^n$ is strong coinitial.
\end{lemma}
{\sc Proof.}
For every $\alpha=(a_1^{p_1}, a_2^{p_2}, \cdots, a_n^{p_n})$ 
objects of the category $\fF T^n/\alpha$ may be considered 
as morphisms
$$
(a_1^{\varepsilon_1}, a_2^{\varepsilon_2}, \cdots, a_n^{\varepsilon_n})
 \rightarrow  
 (a_1^{p_1}, a_2^{p_2}, \cdots, a_n^{p_n}),
$$ 
where $\varepsilon_i\in \{0,1\}$.
The following commutative triangles 
 correspond to it morphisms
$$
  \xymatrix{
	 (a_1^{\varepsilon_1}, a_2^{\varepsilon_2},
         \cdots, a_n^{\varepsilon_n})
	\ar[rr] 
	\ar[dr]
	&& (a_1^{p_1}, a_2^{p_2}, \cdots, a_n^{p_n})\\
	& (a_1^{\varepsilon'_1}, a_2^{\varepsilon'_2}, \cdots, 
		a_n^{\varepsilon'_n}) \ar[ur]	
  }
$$
where $\varepsilon_i\in \{0,1\}$, $\varepsilon'_i\in \{0,1\}$, 
for $1\leq i\leq n$.
We conclude from this that
${\fF T^n}/\alpha \cong (\fF T_1/a_1^{p_1})\times \cdots
 \times(\fF T_n/a_n^{p_n})$ where ${\fF T_i}$ a category consisting 
of twoo the objects
 $1$ and $a_i$
 and two the morphisms $1\stackrel{(1,a_i)}\rightarrow a_i$ and
$1\stackrel{(a_i,1)}\rightarrow a_i$ except the identity morphisms.
Using Lemma \ref{strcof} we get
the isomorphisms $H_q({\fF T^n}/\alpha) \cong H_q({pt})$. 
Consequently, the inclusion $\fF(T^n)\rightarrow \fF\NN^n$ is 
strong coinitial.
\hfill $\Box$

\medskip
\noindent
{\bf Leech homology as a cubical homology.}
 Consider a precubical set $T^n_*$ consisting of sets
$$
T^n_k = \{a_{i_1}a_{i_2}\cdots a_{i_k} : 
1\leq i_1 < i_2 < \cdots < i_k \leq n \}
$$
and maps, for $1\leq s\leq k$,
$$
T^n_{k-1} ~~ { {{\partial^0_s}\atop\longleftarrow}\atop
{\longleftarrow\atop{\partial^1_s}} } ~~ T^n_k ~~,  ~~ 
\partial_s^0(a_{i_1}\cdots
a_{i_k}) = \partial_s^1(a_{i_1}\cdots a_{i_k}) =
a_{i_1}\cdots a_{i_{s-1}}\widehat{a_{i_s}}a_{i_{s+1}} \cdots a_{i_k}
$$

Objects of the category
$\Box_+/T^n_*$ may be considered
as pairs $(k, a)$ where $a\in T^n_k$. Each $a\in T^n_k$ has an unique 
decomposition
 $a_{i_1}\cdots a_{i_k}$ for which 
$1\leq i_1 < \cdots < i_k\leq n$. Thus the objects $(k,a)$ may be identified 
with the elements $a_{i_1}\cdots a_{i_k}\in T^n$. 
Morphisms of  $\Box_+/T^n_*$ are triples 
$(\delta: I^m \rightarrow I^k, a_{j_1}\cdots a_{j_m},
a_{i_1}\cdots a_{i_k})$
satisfying $T^n(\delta)(a_{j_1}\cdots a_{j_m})= a_{i_1}\cdots a_{i_k}$.

It allows us to build a functor 
 $\fS: \Box_+/T^n_* \rightarrow  \fF T^n$ 
as follows. 
Define at objects
$\fS(a_{i_1}\cdots a_{i_k})=a_{i_1}\cdots a_{i_k}$.
Define $\fS$ at morphisms as
$$
 \fS(\delta_s^{k,0}, a_{i_1}\cdots \widehat{a_{i_s}} \cdots a_{i_k}, 
a_{i_1}\cdots a_{i_k}) = 
(a_{i_s},1): a_{i_1}\cdots \widehat{a_{i_s}} \cdots a_{i_k}
\rightarrow 
a_{i_1}\cdots a_{i_k}~,
$$ 
$$
\fS(\delta_s^{k,1}, a_{i_1}\cdots \widehat{a_{i_s}} \cdots a_{i_k}, 
a_{i_1}\cdots a_{i_k}) = 
(1, a_{i_s}): a_{i_1}\cdots \widehat{a_{i_s}} \cdots a_{i_k}
\rightarrow 
a_{i_1}\cdots a_{i_k}~.
$$ 
Every morphism of the category
$\Box_+/T^n_*$ has a decomposition of the form
$(\delta_s^{k,\varepsilon}, a_{i_1}\cdots \widehat{a_{i_s}} \cdots a_{i_k}, 
a_{i_1}\cdots a_{i_k})$, $\varepsilon\in \{0,1\}$. 
There are the commutative diagrams
far all $1\leq s< t\leq k$ and 
$\alpha, \beta \in \{0,1\}$ 
$$
\begin{CD}
a_{i_1}\cdots \widehat{a_{i_s}}\cdots \widehat{a_{i_t}} \cdots a_{i_k}
@>\fS(\delta_s^{\alpha})>>
a_{i_1}\cdots  \widehat{a_{i_t}} \cdots a_{i_k}\\
@V\fS(\delta^\beta_{t-1})VV @VV\fS(\delta^\beta_t) V\\
a_{i_1}\cdots  \widehat{a_{i_s}} \cdots a_{i_k}
@>\fS(\delta^\alpha_s)>>
a_{i_1} \cdots a_{i_k}
\end{CD}
$$
It follows from this that the map  $\fS$ defined above at morphisms 
as
$$
(\delta_s^{k,\varepsilon}, a_{i_1}\cdots \widehat{a_{i_s}} \cdots a_{i_k}, 
a_{i_1}\cdots a_{i_k}), 
$$
has an unique extension by a functor. 
Here  $\fS(\delta)$ is a short denotation 
of the values of $\fS$  at the morphism
$(\delta:
I^m \rightarrow I^k, a_{j_1}\cdots a_{j_m}, a_{i_1}\cdots a_{i_k})$.
It easy to see that the values $\fS(\delta)$ consist of pairs 
of morphisms
$(f,g)$ where $f$ is the value of $\fS$ at the composition 
of morphisms of the form $\delta^0_s$ appeared in the canonical 
decomposition of 
$\delta$ and $g$ is the similar value at the composition of 
$\delta^1_s$.

\begin{lemma}\label{miscpremain}
The functor $\fS: \Box_+/T^n_* \rightarrow  \fF T^n$ is strong 
coinitial.
\end{lemma}
{\sc Proof.} 
Consider the category $\fS/a$. Every morphism $\fS(I^m)\rightarrow a$ 
has the unique decomposition
$$
  \xymatrix{
	\fS(I^m)	\ar[rr] ^{(f,g)}
	\ar[dr]_{\fS(\gamma)}
	&& a\\
	& \fS(I^k) \ar[ur]_{(1,1)}
  }
$$
Consequently  $\fS/a \cong \Box_+/I^k$. The category  
$\Box_+/I^k$ has the terminal object. Hence it is 
connected and acyclic. Thus the functor 
$\fS$ is strong coinitial. 
\hfill $\Box$

We assign to an arbitrary contravariant system 
$F: (\fF\NN^n)^{op} \rightarrow \Ab$  
a homological system 
$\overline{F}= F{\vert}_{{(\fF T^n)^{op}}}\circ \fS^{op}$ 
on the precubical set $T_*^n$.

\begin{theorem}\label{premain}
Let $F: (\fF\NN^n)^{op} \rightarrow \Ab$ be a 
contravariant natural system of abelian groups on the monoid
$\NN^n$. Then for all 
$k\geq 0$ it are hold isomorphisms of 
Leech homology groups and 
cubical homology groups
$$
H_k(\NN^n,F)\cong H_k(T_*^n, \overline{F}).
$$
\end{theorem}
{\sc Proof.} 
The inclusion  $\fF T^n \subset \fF\NN^n$ is strong coinitial 
by Lemma \ref{confn} and hence
$\coLim_n^{(\fF\NN^n)^{op}} F\cong \coLim_n F\vert_{(\fF T^n)^{op}}$.

The functor $\fS$ is strong coinitial by Lemma \ref{miscpremain}.
We have got an isomorphism
$\coLim_k^{(\fF\NN^n)^{op}} F\cong \coLim_k^{(\Box_+/T_*^n)^{op}} 
F{\vert}_{{(\fF T^n)^{op}}}\circ \fS^{op}$. 
The definition of the Leech homology groups 
with the cubical homology groups and Theorem \ref{comcub} 
leads us to the wanted assertion.
\hfill $\Box$

\begin{corollary}
Under the conditions of Theorem \ref{premain}  Leech homology groups
$H_k(\NN^n,F)$ are isomorphic to the 
$k$-th homology groups of the chain complex
$$
0 \leftarrow F(1) \stackrel{d_1}\leftarrow \bigoplus\limits_{1\leq i_1\leq n} F(a_{i_1})
\stackrel{d_2}\leftarrow \bigoplus\limits_{1\leq i_1<i_2\leq n} F(a_{i_1}a_{i_2})
\leftarrow  \cdots
\stackrel{d_n}\leftarrow F(a_1 a_2 \cdots a_n) \leftarrow 0,
$$
where $a_1$, $a_2$, \ldots $a_n$ are generators of the monoid $\NN^n$ and
$$
d_k ((i_1, \cdots, i_k),\varphi) = \sum\limits_{s=1}^k 
(-1)^s((i_1, \cdots, \widehat{i_s}, \cdots, i_k), 
$$
$$
F(a_{i_1}\cdots \widehat{a_{i_s}}\cdots a_{i_k}
\stackrel{(1,a_{i_s})}\rightarrow
 a_{i_1}\cdots  a_{i_k})
 \varphi-F(a_{i_1}\cdots \widehat{a_{i_s}}\cdots
a_{i_k}\stackrel{(a_{i_s},1)}\rightarrow
 a_{i_1}\cdots  a_{i_k})
 \varphi)
$$
\end{corollary}

\begin{corollary}
If for any $1\leq i \leq n$ and $a=a_{i_1}\cdots a_{i_k}$, 
$1\leq i_1<i_2<\cdots <i_k\leq n$ 
it is true  an equality of homomorphisms 
$$
F((1,a_i): a\rightarrow a_i a )= F((a_i,1): a\rightarrow a_i a )
$$ 
for all $i\not\in \{i_1, \ldots, i_k\}$, then   
$H_k(\NN^n,F)= \bigoplus\limits_{1\leq i_1<i_2<\cdots<i_k\leq n} 
F(a_{i_1}a_{i_2}\cdots a_{i_k})$
\end{corollary}
{\sc Proof.} We have $d_k=0$ for all $1\leq k\leq n$. 
This leads to  required isomorphisms.
\hfill $\Box$

\subsection{Leech dimension of a free partially 
commutative monoid}

Now we compute a homological Leech dimension 
of a free partially commutative monoid.

\medskip
\noindent
{\bf Homological dimension of small categories.}
Let $\mC$ be a small category and
$\coLim_n^\mC: \Ab^\mC\rightarrow \Ab$ the left
satellites of the colimit. 
A {\em homological dimension} of $\mC$ 
is defined by
$$
	\hd\mC = \sup\{n\in \NN: \coLim_n^{\mC}\not= 0\}.
$$
Here sup is taken in the totally ordered set
$\{-1\}\cup \NN \cup \{-\infty\}$. For example, $\hd\mC=-1$ 
if and only if $\mC=\emptyset$.

An {\em open subcategory} of the category $\mC$ 
is a full subcategory
$\mD\subseteq \mC$ which contains with each 
its object $d\in \mD$ all the objects
 $c\in \mC$ for which there are morphisms $d\rightarrow c$.
A subcategory  $\mD\subseteq\mC$ is  {\em closed} if  
$\mD^{op}$ is open subcategory of $\mC^{op}$.
\begin{lemma}\label{mitsup}
Let $\mC=\bigcup\limits_{j\in J}\mC_j$ is
an union of some open subcategories
of $\mC$. Then 
$\hd\mC=\sup\limits_{j\in J} \{\hd\mC_j\}$.
\end{lemma}  
{\sc Proof.} For any $a\in \mC$ 
denote by $\mC^a\subseteq \mC$ a full subcategory of $\mC$ which objects  
 $b\in \mC^a$ have morphisms $a\rightarrow b$. It is the smallest open 
subcategory containing $a$.
It follows from Mitchell assertion 
\cite[Corollary 10]{mit1981} that $\hd\mC=\sup\limits_{a\in Ob\mC}\{\hd\mC^a\}$.
Since the subcategories
$\mC_j$ are open, it is true that $\mC^a=\mC^a_j$
for all $a\in \mC_j$. The required equality follows from
$$
\sup\limits_j \{\sup\limits_{a\in Ob\mC_j} \{\hd\mC^a\}\} = 
\sup\limits_{a\in \cup_{j\in J}Ob\mC_j} \{\hd\mC^a\}.
$$
\hfill $\Box$

\medskip
\noindent
{\bf Category of factorizations of a free partially commutative monoid.}
Let $M(E,I)$ be a free partially commutative monoid generated by a set
 $E$ and  relations $ab=ba$ which hold true for all $(a,b)\in I$. 
We say that a subset $S\subseteq E$ consists of 
{\em mutually commuting generators} if $(a,b)\in I$ for all $a,b\in S$.   

Denote by $E_v\subseteq E$, $v\in V$, maximal subsets 
of the mutually commuting generators. Let $M(E_v)\subseteq M(E,I)$ 
be submonoids generated by $E_v$.
Denote by  $E_v\subseteq E$, $v\in V$, the maximal subsets 
consisting of its mutually independent generators.
Let $M(E_v)\subseteq M(E,I)$ be denote the monoid generated by
 $E_v$.

The union $\bigcup\limits_{v\in V} M(E_v)\subseteq  M(E,I)$ is not a subcategory.
It makes difficulties for the studying homology of $M(E,I)$. 
The idea to considering the inclusion of the 
corresponding categories of factorizations  saves the situation:
 The full subcategory of $\fF M(E,I)$ with the set of objects 
$M(E_v)$ equals  $\fF M(E_v)$.
Moreover, the subcategories $\fF M(E_v)\subseteq \fF M(E,I)$ are closed.

First we study the comma category $\fF M(E,I)/\alpha$ for any element  
$\alpha\in M(E,I)$. 

\begin{proposition}\label{comp}
Let $M(E,I)$ be a free partially commutative monoid and $\alpha \in M(E,I)$ an element. 
Then the category $\fF M(E,I)/\alpha$ is a partially ordered set.
\end{proposition}
{\sc Proof.} By  \cite[Corollary 2]{die1997} if $E$ is finite, then
the monoid 
 $M(E,I)$ is cancellative. Let  $E$ be infinite. Suppose that
$y w_1 x = y w_2 x$. There is a finite subset $E'\subseteq E$ consisting of symbols 
 which are multipliers of the elements $x, y, w_1, w_2 \in M(E,I)$.
Since $E'$ is finite, the submonoid $M(E', I\cap(E'\times E'))$ is cancellative. It follows that
 $w_1=w_2$. Hence,  $M(E,I)$ is cancellative. 
Assign to every  $x\in M(E,I)$ a length 
$|x|$  
of a word which presents $x$. It is called a {\em length} of $x$.  
It easy to see that $|xy|=|x|+|y|$ for all $x,y \in M(E,I)$. 
Therefore,  $M(E,I)$ does not contain nonidentity retractions. 
By Lemma \ref{commapart} we conclude that $\fF M(E,I)/\alpha$ is a partially ordered set.
\hfill $\Box$

\begin{example}
\em
Let $M(E,I)$ be a free partially commutative monoid and 
$a\in E$ a generator.
Describe the partially ordered set $\fF M(E,I)/a$. It consists of the elements $(1,1,a)$, 
$(1,a,1)$, and $(a,1,1)$. Its morphisms are defined by the diagram
$$
  
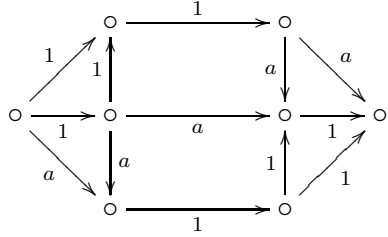
\xymatrix{
	& \circ	\ar[rr] ^{1}  & & \circ \ar[dr] ^{a} \ar[d]_{a}\\
   \circ \ar[dr]_a \ar[ur]^1 \ar[r]_1 & \circ \ar[d]^a \ar[u]^1
   \ar[rr]_a && \circ\ar[r]_1 & \circ \\
	& \circ \ar[rr]_{1} & & \ar[u]^1 \circ \ar[ur]_1 
  }
$$
The morphisms $(1,1,a) \stackrel{(1,a)}\rightarrow (1,a,1)$ and 
$(a,1,1) \stackrel{(a,1)}\rightarrow (1,a,1)$ leads us to the relation 
$(1,1,a) < (1,a,1) > (a,1,1)$  in the partially ordered set
$\fF M(E,I)/a$. 
\end{example}

Elements  $x, y \in M(E,I)$ are said to be {\em commuting} if $x y= y x$.
We notice attention that generators $a$ and $b$ are commuting if and only if $(a,b)\in I$ or $a=b$.

The category $\fF M(E,I)$ contains the full subcategory 
$\bigcup\limits_{v\in V}\fF M(E_v)$ whose objects are products 
of the pairwise commuting generators. 

\begin{theorem}\label{coinit}
The inclusion $\bigcup\limits_{v\in V}\fF M(E_v)\subseteq \fF M(E,I)$ 
is strong coinitial.
\end{theorem}
{\sc Proof.} 
By Lemma \ref{factop} the category $\fF M(E,I)$ is isomorphic to 
$\fF (M(E,I)^{op})$. This isomorphism carries the subcategory 
$\bigcup\limits_{v\in V}\fF M(E_v)$ to   
$\bigcup\limits_{v\in V}\fF (M(E_v)^{op})$. Thus, it is enough to prove that
the inclusion 
$\bigcup\limits_{v\in V}\fF (M(E_v)^{op})\subseteq \fF (M(E,I)^{op})$ 
is strong coinitial.
It is equivalent to the assertion that for each
$\alpha \in M(E,I)$
$$
	H_q(\bigcup\limits_{v\in V}\fF (M(E_v)^{op})/\alpha)\cong H_q(pt), \mbox{for all} ~ q\in \NN~.
$$
We prove this assertion by induction on the length
$\vert\alpha\vert$.
By Proposition \ref{comp} the comma category $\fF (M(E,I)^{op})/\alpha$ 
is a partially ordered set.
By Corollary \ref{factopcor} its elements may be given as  
triples  $(x,\beta,y)$ of elements of $M(E,I)$ 
satisfying $x\circ\beta\circ y = \alpha$.
The relation $(x,\beta,y)\leq (x',\beta',y')$ is defined by existence of pairs 
 $(v,w)$ making the following diagram to be commutative
$$
	\xymatrix{
		& \circ \ar[ld]_x \ar[dd]_v & \circ \ar[l]^\beta 	\\
		\circ &&& \circ \ar[lu]_y \ar[ld]^{y'}\\
		& \circ \ar[lu]^{x'} & \circ \ar[l]_{\beta'} \ar[uu]_w
	}
$$
The subcategory $(\bigcup\limits_{v\in V}\fF (M(E_v)^{op}))/\alpha \cong 
\bigcup\limits_{v\in V} (\fF (M(E_v)^{op})/\alpha)$ is its closed 
subset consisting of the triples $(x,\beta,y)$ for which $\beta \in \bigcup\limits_{v\in V} M(E_v)$.

Let $\Phi(\alpha)=\bigcup\limits_{v\in V} (\fF (M(E_v)^{op})/\alpha)$. 
We will prove by induction on $\vert\alpha\vert$ that the partially ordered set 
$\Phi(\alpha)$ is acyclic. 
We use that  $\Phi(\alpha)$
is union of closed subsets which are isomorphic to  
$\Phi(\beta)$ with $\vert\beta\vert < \vert\alpha\vert$.

If $\vert\alpha\vert=0$, then $\alpha=1$, hence the partially ordered set   
$\Phi(\alpha)$ consists of unique element
$(1,1,1)$. So, the assertion is true in this case.

Suppose that $\forall \beta\in M(E,I)$ with $\vert\beta\vert < n$ 
there are isomorphisms $H_q(\Phi(\beta)) \cong 
H_q(pt)$ for all $q\geq 0$. Prove it for $\alpha$ with $\vert\alpha\vert=n$.

If $\alpha=f\circ g$, then  we say that $f$ is a {\em left divisor} of  $\alpha\in M(E,I)$ and 
$g$ is its {\em right divisor}. In this case,
there are  injections of the partially ordered sets 
$$
	\fF (M(E,I)^{op})/g \stackrel{f_*}\longrightarrow \fF (M(E,I)^{op})/\alpha
		\stackrel{g^*}\longleftarrow \fF (M(E,I)^{op})/f,
$$
where $f_*(x,\beta,y)=
(f\circ x, \beta,y)$, 
$g^*(x,\beta,y)=
(x, \beta, y\circ g)$.
Images of the maps 
$$
\Phi(g) \stackrel{f_*}\rightarrow
\Phi(\alpha) \stackrel{g^*}\leftarrow \Phi(f)
$$
are closed in
$\Phi(\alpha)$. 
In general, if $f\circ\gamma\circ g = \alpha$ for some  $f,g,\gamma \in M(E,I)$, then
the image of the inclusion
$$
f_* g^* = g^* f_*: 
\Phi(\gamma)
\longrightarrow
\Phi(\alpha)
$$
is closed in  
$\Phi(\alpha)
$.

Let $\{a_1, \cdots, a_m\}\subseteq E$ be the set of generators 
which are left divisors of $\alpha$ and $\{b_1, \cdots , b_n\}\subseteq E$
 its right divisors.
It is easy to see that the generators  $a_1, \cdots, a_m$ are mutually independent.
Similarly $b_1, \cdots , b_n$ are mutually independent.
Denote by 
$$
L_i(\alpha)=\{(x,\beta,y)\in \Phi(\alpha)
: 
 x \mbox{  has the left divisor } a_i\}
$$
$$
R_j(\alpha)=\{(x,\beta,y)\in \Phi(\alpha)
: 
 y \mbox{  has the right divisor } b_j\}
$$
for every $1\leq i\leq m$ and $1\leq j\leq n$.
Let $L(\alpha)=\bigcup\limits_{i=1}^m L_i(\alpha)$ and  
 $R(\alpha)=\bigcup\limits_{j=1}^n R_j(\alpha)$.
It is evidently that $L(\alpha)= \{(x,\beta,y): \beta\in \bigcup\limits_{v\in V} M(E_v), x\not=1 \}$ and
$R(\alpha)= \{(x,\beta,y): \beta\in \bigcup\limits_{v\in V} M(E_v), y\not=1 \}$.

Now consider the following two cases: $\alpha\in \bigcup\limits_{v\in V} M(E_v)$ and 
$\alpha \not\in \bigcup\limits_{v\in V} M(E_v)$. 
If $\alpha\in \bigcup\limits_{v\in V} M(E_v)$, then 
$\Phi(\alpha)$ contains the largest element
$(1, \alpha, 1)$ and consequently 
 $H_q(\Phi(\alpha))\cong H_q(pt)$ for all $q\geq 0$.

If 
$\alpha \not\in \bigcup\limits_{v\in V} M(E_v)$, then $x\not=1$ or $y\not=1$ holds for each 
$(x, \beta, y) \in \Phi(\alpha).
$
Therefore 
 $(x,\beta,y)\in L(\alpha)\cup R(\alpha)$. It follows that  
$\Phi(\alpha)= L(\alpha)\cup R(\alpha)$.

Now we will prove the acyclic property of $L(\alpha)$, $R(\alpha)$, and $L(\alpha)\cap R(\alpha)$.
Then, by Lemma \ref{covlem}, we will obtain the acyclic property of $L(\alpha)\cup R(\alpha)$.
If  $\alpha=f\circ g$, then denote   
$g$ by $f^{-1}\alpha$ and $f$ by $\alpha g^{-1}$.

\begin{itemize}
\item $L(\alpha)=\bigcup\limits_{i=1}^m L_i(\alpha)$, 
$L_{i_1}(\alpha)\cap\cdots \cap L_{i_k}(\alpha)=(a_{i_1}\cdots a_{i_k})_*\Phi((a_{i_1}\cdots a_{i_k})^{-1}\alpha)$.
Since the map $(a_{i_1}\cdots a_{i_k})_*$ makes an isomorphism of the partially ordered set
 $\Phi((a_{i_1}\cdots a_{i_k})^{-1}\alpha)$ with a closed subset, the intersections 
$L_{i_1}(\alpha)\cap\cdots \cap L_{i_k}(\alpha)$ are acyclic. By Corollary \ref{covcor} we obtain 
the acyclic property of
 $L(\alpha)$.
\item  The acyclic property of $R(\alpha)=\bigcup\limits_{j=1}^n R_j(\alpha)$ is proved similarly.
\item $L(\alpha)\cap R(\alpha)= \bigcup\limits_{i=1}^m (L_i(\alpha)\cap R(\alpha))$, \quad
$L_{i_1}(\alpha)\cap\cdots \cap L_{i_k}(\alpha)\cap R(\alpha)= 
\{(a_{i_1}\cdots a_{i_k}x, \beta, y): \beta\in \bigcup\limits_{v\in V}M(E_v), y\not=1, 
a_{i_1}\cdots a_{i_k}x \beta y = \alpha\}$. It follows that 
$$
L_{i_1}(\alpha)\cap\cdots \cap L_{i_k}(\alpha)\cap R(\alpha)=
(a_{i_1}\cdots a_{i_k})_*R((a_{i_1}\cdots a_{i_k})^{-1}\alpha)
$$
is isomorphic to the partially ordered set $R((a_{i_1}\cdots a_{i_k})^{-1}\alpha)$ whose acyclic property 
is proved in the second item.
\end{itemize}

It follows from Lemma \ref{covlem} that $\Phi(\alpha)$ is acyclic.
\hfill $\Box$

\begin{definition}
Let $M$ be a monoid. Its a {\em homological Leech dimension} is defined by
 $\Ld M = \hd (\fF M)^{op}$.
\end{definition}

\begin{corollary}\label{hyptes}
Let $M(E,I)$ be a free partially commutative monoid 
in which the maximal number of mutually commuting generators 
equals  $n$. 
Then $\Ld M(E,I)=n$. 
\end{corollary}
{\sc Proof.} 
Since $\bigcup\limits_{v\in V} \fF M(E_v) \subseteq \fF M(E,I)$ is 
the closed subcategory, the following inequality 
follows from the assertion
\cite[Corollary 10]{mit1981} 
$$
\hd (\bigcup\limits_{v\in V} \fF M(E_v))^{op} \leq \hd (\fF M(E,I))^{op}.
$$
The inclusion 
$\bigcup\limits_{v\in V} \fF M(E_v) \subseteq \fF M(E,I)$
is strong coinitial by Theorem \ref{coinit}. Thus 
for any functor $F: (\fF M(E,I))^{op}\rightarrow \Ab$ it follows from 
$\coLim_n^{(\fF M(E,I))^{op}} F\not=0$ that 
$\coLim_n F\vert_{(\bigcup\limits_{v\in V} \fF M(E_v))^{op}} \not=0$.
From here 
$$
\hd (\bigcup\limits_{v\in V} \fF M(E_v))^{op} \geq \hd (\fF M(E,I))^{op}.
$$
Consequently, 
$\hd (\bigcup\limits_{v\in V} \fF M(E_v))^{op} = \hd (\fF M(E,I))^{op}$.
From this, using Lemma \ref{mitsup} we get the equality
$\hd (\fF M(E,I))^{op}=\sup\limits_{v\in V}\hd(\fF M(E_v)^{op})$.
So
 $\Ld M(E,I)= \sup\limits_{v\in V}\Ld M(E_v)$.
Since the homological Leech dimension 
of the free commutative monoid generated by 
 $n$ lements equals  $n$, we conclude 
that $\Ld M(E,I)=n$ where $n$ is the maximal number of 
generators of  $M(E_v)$.
\hfill $\Box$

\subsection{Leech homology groups of free partially commutative 
monoids}

\medskip
\noindent
{\bf A coinitial subcategory of the factorization category.}
Let $E$ be a set totally ordered by a relation  $<$. 
We prove the following auxiliary assertion.
\begin{proposition}
Let  $T_v\subset M(E_v)$, for every $v\in V$, be a set  
consisting of products
$a_1 a_2 \cdots a_n$ such that $a_1 < a_2 < \cdots < a_n$ and  
$a_j\in E_v$ for all $1\leq j\leq n$. Here  $n$ is 
a finite number such that 
$n \leq\vert E_v \vert$.  The product equals $1\in T_v$ for  $n=0$.
Then the inclusion  $\bigcup\limits_{v\in V} \fF T_v \subset \fF M(E,I)$ 
is strong coinitial.
\end{proposition}
{\sc Proof.} A composition of strong coinitial 
functors is strong coinitial by Corollary \ref{compos}.
We have proved 
that the $\bigcup\limits_{v\in V} \fF M(E_v) \subseteq \fF M(E,I)$ 
is strong coinitial. 
Thus it is enough to show that the inclusion
$\bigcup\limits_{v\in V} \fF T_v \subset\bigcup\limits_{v\in V} \fF M(E_v)$
 is strong coinitial.
But for each $\alpha \in \bigcup\limits_{v\in V} \fF M(E_v)$ there exists 
$w\in V$ such that $\alpha \in \fF M(E_w)$. All factors of $\alpha$
belong to $M(E_w)$. Consequently
$\bigcup\limits_{v\in V} \fF T_v/\alpha = \fF T_w/\alpha$.

The inclusion  
$\fF T_w\subset \fF M(E_w)$ is strong coinitial by Lemma 
\ref{confn}. So,  
$H_q(\bigcup\limits_{v\in V} \fF T_v/\alpha)\cong H_q(pt)$.
\hfill $\Box$

\medskip
\noindent
{\bf Main Theorem.}
Let $E$
be a set and $I\subseteq E\times E$ an irreflexive symmetric relation.
We build a precubical set 
$T(E,I)$ depending on some total ordering relation $\leq$ on $E$.
For every integer  $n>0$ 
define  $T(E,I)_n$ as a set of words $a_1\cdots a_n$ consisting 
of mutually commuting generators $a_1<\cdots < a_n$ where 
$a_1, \cdots, a_n \in E$ 
$$
T(E,I)_n = \{a_1 \cdots a_n: (a_1 < \cdots < a_n) \& 
( 1\leq i<j\leq n \Rightarrow (a_i, a_j)\in I) \}.
$$
(For $n=0$ the set  $T(E,I)_0$ consists of the empty word $1$.)
The maps $\partial^{n,\varepsilon}_i: T(E,I)_n\rightarrow T(E,I)_{n-1}$, 
for $1\leq i\leq n$, act as
$$
\partial^{n,0}_i(a_1 \cdots  a_n) = \partial^{n,1}_i(a_1 \cdots a_n) = 
a_1 \cdots  \widehat{a_i}  \cdots a_n~.
$$

It easy to see that $T(E,I)$ is union of precubical sets
$(T_v)_*$ for which  
$$
(T_v)_n = \{a_1\cdots a_n \in T(E,I)_n~:~ ( 1\leq i\leq n \Rightarrow
a_i\in E_v)\},
$$
where $E_v\subseteq E$ are the maximal subsets 
of mutually commuting generators 
of the monoid $M(E,I)$.
The maps
$(T_v)_n \stackrel{\partial^{n,\varepsilon}_i}\rightarrow (T_v)_{n-1}$,
act for $1\leq i \leq n$ and $\varepsilon \in \{0,1\}$
 by
$$
\partial^{n,0}_s(a_1\cdots a_n) =  \partial^{n,1}_s(a_1\cdots a_n) =
a_1\cdots \widehat{a}_s \cdots a_n
$$
Let $\fS: \Box_+/\bigcup\limits_{v\in V} (T_v)_* \rightarrow 
\bigcup\limits_{v\in V} \fF T_v$ be a functor which 
assign to every singular cube 
 $a_1\cdots a_n \in \bigcup\limits_{v\in V} (T_v)_n$
the object $a_1\cdots a_n$. The action of functor $\fS$ at morphisms is 
similar with the action in the case of free commutative monoid.

For any functor
$F: (\fF M(E,I))^{op} \rightarrow \Ab$ 
denote by $\overline{F}$
a homological system on $T(E,I)$ defined as the composition
$$
(\Box_+/T(E,I))^{op} \stackrel{\fS^{op}}\longrightarrow 
{\bigcup\limits_{v\in V}(\fF T_v)^{op}} \subset (\fF M(E,I))^{op}
\stackrel{F}\rightarrow  \Ab
$$

\begin{definition}
A free partially commutative monoid  $M(E,I)$ is called 
{\em locally bounded} if every set of mutually commuting generators 
is finite.  
\end{definition}

\begin{theorem}\label{leechiscubic}
Let $M(E,I)$ be a locally bounded free partially commutative monoid an 
$F: (\fF M(E,I))^{op} \rightarrow \Ab$ a functor.
Then for all $n\geq 0$ there are isomorphisms 
of the Leech and cubical homology groups
$$
H_n(M(E,I),F) \cong H_n(T(E,I),\overline{F}).
$$
\end{theorem}
{\sc Proof.}  
The functor
$\fS$ carries the subcategory $\Box_+/(T_v)_*\subseteq \Box_+/T(E,I)$
into  $\fF T_v$ and hence defines a functor which 
we denote by
$\fS_v: \Box_+/ (T_v)_* \rightarrow \fF T_v$.
For an arbitrary
$\alpha \in \bigcup\limits_{v\in V} \fF T_v$
there exists $w\in V$ such that
$\alpha \in \fF T_w$. Consequently there is an isomorphism
$\fS/\alpha \cong \fS_w/\alpha$.
It follows from Lemma \ref{miscpremain} 
for the free finite generated commutative monoid 
$M(E_v)$ that the functor $\fS_v: \Box_+/ (T_v)_* \rightarrow  \fF T_v$ 
is strong coinitial. So the category  
$\fS_v/\alpha$ is connected and acyclic. Therefore the category
$\fS/\alpha$ is connected and acyclic. Thus  $\fS$ is strong coinitial.
Using that the inclusion 
$\bigcup\limits_{v\in V} \fF T_v \subseteq \fF M(E,I)$ is strong 
coinitial we get a required isomorphisms
$$
\coLim_n^{(\fF M(E,I))^{op}} F \cong 
\coLim_n^{(\Box_+/\bigcup\limits_{v\in V} 
(T_v)_*)^{op}}~\left(F\vert_{\bigcup\limits_{v\in V}(\fF T_v)^{op}}
~\circ\fS^{op}\right)
$$
\hfill $\Box$

From this we obtain the main result of the paper

\begin{theorem}\label{main}
Let $M(E,I)$ be a locally bounded free partially commutative monoid and  
$F: (\fF M(E,I))^{op} \rightarrow \Ab$ a functor. 
Then groups $\coLim_n^{(\fF M(E,I))^{op}}F$
are isomorphic to $n$-th homology groups of a chain complex
$$
0 \leftarrow F(1) \stackrel{d_1}\leftarrow
\bigoplus\limits_{a_1\in E} F(a_1)
\stackrel{d_2}\leftarrow
 \bigoplus\limits_{{a_1<a_2}\atop{(a_1,a_2)\in I}} F(a_1a_2)
\leftarrow  \cdots ~~~~~~~~~~~~~~~~~~~~~~~~~~~~~~~~
$$
$$
~~~\cdots ~~~ \leftarrow
\bigoplus\limits_{{a_1<a_2<\cdots<a_{n-1}}\atop{(a_i,a_j)\in I} } 
F(a_1 a_2\cdots a_{n-1}) 
 \stackrel{d_n}\longleftarrow 
\bigoplus\limits_{{a_1<a_2<\cdots<a_n}\atop{(a_i,a_j)\in I} } 
F(a_1 a_2\cdots a_n) \leftarrow \cdots~,
$$
$n$-th member of which for  $n\geq 0$ equals a direct sum 
of abelian groups 
$F(a_1\cdots a_n)$ 
taken by all sequences  
$a_1<a_2<\cdots<a_n$ of mutually commuting generators 
$a_1, a_2, \cdots, a_n \in E$
with the differentials
$$
d_n((a_1, \cdots, a_n),\varphi) = 
\sum\limits_{s=1}^n (-1)^s ((a_1, \cdots, \widehat{a_s} , \cdots, a_n),
F(1,a_s)\varphi-F(a_s,1)\varphi)
$$
where $F(1,a_s)$ shortly denote values $F$ at morphisms 
$(1,a_s):a_1\cdots\widehat{a_s}\cdots a_n \rightarrow a_1\cdots a_n$, 
and $F(a_s,1)$ at 
$(a_s,1):a_1\cdots\widehat{a_s}\cdots a_n \rightarrow a_1\cdots a_n$.
\end{theorem}
{\sc Proof.}
It follows from Theorem \ref{leechiscubic} that 
$\coLim_n^{(\fF M(E,I))^{op}}F \cong 
\coLim_n^{(\Box_+/T(E,I))^{op}}~\overline{F}$.
Theorem \ref{comcub} leads to the isomorphisms
$$
\coLim_n^{\fF M(E,I)^{op}}F \cong 
H_n(C_*(T(E,I),\overline{F})).
$$
\hfill $\Box$

\subsection{Hochschild homology}

\begin{lemma}
Let $S: \mC\rightarrow \mD$ be a functor between small categories.
Then there exists a isomorphism which is natural in  
$F\in \Ab^{\mC}$ and $G\in \Ab^{\mD^{op}}$
$$
	\Lan^{S}F \otimes G \cong F\otimes (G\circ S).
$$
\end{lemma}
{\sc Proof.} Recall that  
$\Hom(F,A): \mC^{op}\rightarrow \Ab$ for an abelian group $A$
denotes a functor having values 
 $\Ab(F(c),A)$ at objects $c\in \mC$ at $\Ab(F(\alpha),1_A)$ at 
morphisms $\alpha \in Mor(\mC)$. Since  $\Lan^{S}$ is left adjoint to 
$(-)\circ S$, there are isomorphisms which are 
natural in  
$F\in \Ab^{\mC}$, $G\in \Ab^{\mD^{op}}$, and $A\in \Ab$ 
$$
  \Ab^{\mD}(\Lan^{S}F, \Hom(G,A)) 
  \cong \Ab^{\mC}(F, \Hom (G, A)\circ S) 
\cong \Ab^{\mC}(F, \Hom (G\circ S^{op}, A)).
$$
It follows from isomorphisms  
$\Ab(T\otimes G, A)\cong \Ab^{\mD}(T, \Hom(G,A))$ 
which are natural in  $T\in \Ab^{\mD}$ 
that there exists the isomorphism  
$$
\Ab(\Lan^{S}F\otimes G, A) \stackrel{\cong}\rightarrow 
 \Ab(F\otimes (G\circ {S^{op}}), A)
$$
which is natural in three arguments.
Since the Yoneda inclusion $T \mapsto \Ab(T,-)$ is full, there exists 
an isomorphism 
$F\otimes (G\circ S) \stackrel{\cong}\rightarrow \Lan^{S}F \otimes G$.
\hfill $\Box$

Let $(s,t): \fF\mC \rightarrow \mC^{op}\times \mC$ be a functor 
assigning to every object
$a\stackrel{\alpha}\rightarrow b$ the pair 
$(a,b)\in Ob(\mC^{op}\times \mC)$ and to every morphism  
$\alpha \stackrel{(f,g)}\rightarrow \beta$ the pair  
$(f,g)\in Mor(\mC^{op}\times \mC)$.

Let $\mC$ be a small category. 
We define its {\em Hochschild--Mitchell homology 
groups with coefficients in a $\mC^{op}$--bimodule 
 $B: (\mC^{op}\times \mC)^{op} \rightarrow \Ab$} 
as 
$\Tor_n(L\mC,B)$ where $L\mC: \mC^{op}\times \mC\rightarrow \Ab$ 
is the composition of functors
$\mC^{op}\times \mC
\stackrel{\mC(-,=)}\rightarrow \Set \stackrel{L}\rightarrow \Ab$.
\begin{lemma}
$\Tor_n(L\mC,B)\cong \coLim_n^{\fF{\mC}^{op}} B\circ{(s,t)^{op}}$
\end{lemma}
{\sc Proof.} Since connected components of 
the category  $(s,t)/(a,b)$ 
have terminals objects by \cite[Lemma 2.4]{X1997}, the application of 
$\Lan^{(s,t)}$ to a projective resulution 
 $P_*\rightarrow \Delta \ZZ$ gives a projective resolution 
of $\Lan^{(s,t)}\Delta\ZZ=L\mC$. From the isomorphisms 
$\Lan^{(s,t)}P_n\otimes B \cong P_n\otimes B\circ (s,t)^{op}$
we obtain the isomorphisms
$$
\Tor_n(\Lan^{(s,t)}\Delta\ZZ, B) \cong \Tor_n(\Delta\ZZ, B\circ (s,t)^{op})
$$
The required assertion follows from  
$\Tor_n(\Delta\ZZ, B\circ (s,t)^{op})\cong 
\coLim_n^{\fF{\mC}^{op}} B\circ{(s,t)^{op}}$.
\hfill $\Box$

Let $M$ be a monoid considered as a category with one 
object. Then its Hochschild--Mitchell homology 
groups are isomorphic to Hochschild homology groups  
$H_n(\ZZ[M],B)$ of the monoid ring  $\ZZ[M]$ with coefficients 
in  
$\ZZ[M]\otimes\ZZ[M]^{op}$--module $B$.

\begin{corollary}
Let $M(E,I)$ be a locally bounded free 
partially commutative monoid
and 
$B$ be a $\ZZ[M(E,I)]\otimes\ZZ[M(E,I)]^{op}$--module.
Then the Hochschild homology groups 
 $H_n(\ZZ[M(E,I)],B)$ are isomorphic to homology groups of the chain complex
consicting of the abelian groups $C_0=B$ and  
$C_n=\bigoplus\limits_{(a_1, \cdots, a_n)}B$ for $n>0$ where the sums 
are taken by all sequences of mutually 
commuting generators  
$a_1< \cdots <a_n$,  $a_1, \cdots, a_n\in E$.
Differentials are defined by the following values, for $\varphi\in B$,
$$
d_n(a_1, \cdots, a_n, \varphi)= 
\sum_{i=1}^n (-1)^i (a_1, \cdots, \widehat{a_i}, \cdots , a_n, 
a_i\cdot\varphi - \varphi\cdot a_i ).
$$
Here the dot denotes both a left and a right action of a monoid on $B$. 
\end{corollary}
(These actions may be written as $a\cdot\varphi\cdot b = B(a,b)(\varphi)$
by the bifunctor $B(-,=)$.)

\subsection{Homology with coefficients in a right module}

 Let $\mC$ be an arbitrary small category. Consider a functor
$\fF\mC \stackrel{t}\rightarrow \mC$ assigning to every
arrow
$a\rightarrow b$ the end $b$ and to every morphism  
$\alpha \stackrel{(f,g)}\rightarrow \beta$ the morphism  
$g: t(\alpha)\rightarrow t(\beta)$. The following assertion 
is implicitly contained in  \cite{X1997}. 

\begin{lemma}\label{hash}
The functor $\fF\mC \stackrel{t}\rightarrow \mC$ is strong 
coinitial.
\end{lemma}
{\sc Proof.} For each object $c\in \mC$ the category $t/c$ 
contains a full subcategory $\mD$ consisting of all objects $(\beta,1_c)$.
This subcategory is isomorphic to
$(\mC/c)^{op}$. For any object
$(\alpha,x)$ it is defined a morphism from   
$(\alpha,x)$ to  $(x\circ\alpha,1_c) \in Ob(\mD)$ given by the 
diagram

$$
 \xymatrix{
	a \ar[rr]^{\alpha} & & b \ar[rrd]^{x} \ar[dd]^x \\
	& & & & c \\
	a \ar[uu]_{1_a} \ar[rr]_{x\circ\alpha} & & c \ar[urr]_{1_c}
 }
$$
The morphism $(\alpha,x)\rightarrow (x\circ\alpha,1_c)$ is 
universal in the sense 
that for each morphism 
$(\alpha,x)\rightarrow (\beta,1_c)$ 
there exists an unique morphism 
$(x\circ\alpha,1_c)\rightarrow (\beta,1_c)$ 
making commutative the following diagram

$$
 \xymatrix{
	(\alpha,x)  \ar[rr] \ar[rd] 
		& & (x\circ\alpha, 1_c) \ar@{-->}[ld]_{\exists!} \\
		& 
	(\beta,1_c)   
 }
$$
Consequently, there is a functor $\sigma: t/c\rightarrow (\mC/c)^{op}$
which is left adjoint  to a full inclusion
$(\mC/c)^{op}\subseteq t/c$.
Hence the inclusion $\mC/c\subseteq (t/c)^{op}$ 
is strong coinitial. Therefore
$\coLim_n^{t/c}(\Delta\ZZ)\circ \sigma
\cong \coLim_n^{(\mC/c)^{op}}(\Delta\ZZ)$.
So
$H_q(t/c)=0$ for $q>0$ and $H_0(t/c)\cong \ZZ$.
\hfill $\Box$

\begin{corollary}\label{hyptes1}
Let $M(E,I)$ be a free partially commutative monoid
 which does not contain   
$>~n$ mutually commuting generators. 
Then  $\hd M(E,I)\leq n$.
\end{corollary}
{\sc Proof.} It follows from Lemma \ref{hash} for $\mC= M(E,I)^{op}$ that 
the functor $t: \fF (M(E,I)^{op}) \rightarrow M(E,I)^{op}$ 
is strong coinitial.
The categories $\fF (M(E,I)^{op})$ and $\fF M(E,I)$ are isomorphic. 
Therefore 
$\hd (\fF M(E,I))^{op}\geq \hd M(E,I)$. By Corollary \ref{hyptes} we have 
 $\hd M(E,I)\leq n$.
\hfill $\Box$

\begin{corollary}\label{hrm}
Let $M(E,I)$ be a locally bounded free partially commutative 
monoid and 
G a right $M(E,I)$--module. Then the homology groups  $H_n(M(E,I)^{op},G)$
are isomorphic to the homology groups of the chain complex
$$
0 \leftarrow G \stackrel{d_1}\leftarrow \bigoplus\limits_{a_1\in E} G
\stackrel{d_2}\leftarrow \bigoplus\limits_{{a_1<a_2}\atop{(a_1,a_2)\in I}} G
\leftarrow  \cdots 
$$
$$
~~~~~~~~~~~~~~~~~~\cdots ~~~\leftarrow  
\bigoplus\limits_{{a_1<a_2<\cdots<a_{n-1}}\atop{(a_i,a_j)\in I} } 
G
\stackrel{d_n}\longleftarrow 
\bigoplus\limits_{{a_1<a_2<\cdots<a_n}\atop{(a_i,a_j)\in I} } 
G \leftarrow \cdots~,
$$
the $n$-th member of which for $n\geq 0$ is equal to 
the direct sum of copies of the group
$G$ 
by all sequences of mutually commuting generators
 $a_1<a_2<\cdots<a_n$ belonging to $E$
and the differential is defined by
$$
d_n((a_1, \cdots, a_n),g) = 
\sum\limits_{s=1}^n (-1)^s ((a_1, \cdots, \widehat{a_s} , \cdots, a_n),
ga_s-g)
$$
where $ga_s\in G$ denotes a member 
obtained by right action of   $a_s$ on   
$g\in G$.
\end{corollary}
{\sc Proof.} By Lemma \ref{hash}  
$t: \fF M(E,I)\rightarrow M(E,I)$ is strong coinitial. Thus the groups
$H_n(M(E,I)^{op},G)$, which  equal 
$\coLim_n^{ M(E,I)^{op}}G$ by definition, are isomorphic 
to
$\coLim_n^{ \fF M(E,I)^{op}}G\circ t^{op}$. Substituting  $t(1,a_s)=a_s$,
$t(a_s,1)=1$, $G\circ t^{op}(1,a_s)(g)=ga_s$, $G\circ t^{op}(a_s,1)(g)=1$, 
we obtain from Theorem \ref{main} the required chain complex.
\hfill $\Box$

\section{Homology of right $M(E,I)$--sets with coefficients in diagrams}

\subsection{Homology groups of right $M(E,I)$--sets} 

\begin{theorem}\label{main2}
Let $M(E,I)$ be a locally bounded free partially commutative monoid, 
$X: M(E,I)^{op}\rightarrow \Set$ a right $M(E,I)$--set, 
$F: (M(E,I)/X)^{op} \rightarrow \Ab$ a functor.
Then the homology groups  $\coLim_n^{(M(E,I)/X)^{op}}F$ are isomorphic to the homology groups 
of the chain complex
$$
0 \leftarrow 
\bigoplus\limits_{x\in X} F(x) \stackrel{d_1}\leftarrow 
\bigoplus\limits_{(x,a_1)} F(x)
\stackrel{d_2}\leftarrow \bigoplus\limits_{{(x, a_1 a_2), a_1<a_2}
\atop{(a_1,a_2)\in I}} F(x)
\leftarrow  \cdots 
$$
$$
~~~~~~~\cdots ~~~\leftarrow  
\bigoplus\limits_{(x, a_1\cdots a_{n-1}), 
{a_1<a_2<\cdots<a_{n-1}}\atop{(a_i,a_j)\in I} } 
F(x)
\stackrel{d_n}\longleftarrow 
\bigoplus\limits_{{(x, a_1\cdots a_n),
 a_1<a_2<\cdots<a_n}\atop{(a_i,a_j)\in I} }
F(x) \leftarrow \cdots~,
$$
with  
$$
d_n(x,a_1\cdots a_n,g) = 
$$
$$
\sum_{s=1}^n(-1)^s \left(
(xa_s, a_1\cdots \widehat{a_s}\cdots a_n, 
F(x\stackrel{a_s}\rightarrow x\cdot a_s)g)-
(x, a_1\cdots \widehat{a_s}\cdots a_n, g) \right)
$$ 
\end{theorem}
{\sc Proof.} In a similar manner  with the proof of
\cite[Theorem 5.3]{X20042} consider an arbitrary 
monoid  $M$ with functors   
 $X: M^{op}\rightarrow \Set$ and  $F: (M/X)^{op}\rightarrow \Ab$.
Let $S: (M/X)^{op}\rightarrow M^{op}$ be a (dual of forgetful) functor  
which acts as $S(x\stackrel{\mu}\rightarrow y)= \mu$. For 
any functor $F: (M/X)^{op}\rightarrow \Ab$ there is the 
left Kan extension $\Lan^S F$. Each of the connected components 
of $S/M$ where  $M$ denotes the unique object of the monoid 
considered as a small category has a terminal object of the form
$(x,1)$ for some $x\in X$. It follows from here that 
 $\Lan^S F$ is a right 
 $M$--module $\bigoplus_{x\in X}F(x)$ with the action 
$(x,g)\cdot\mu=(x\mu, F(x\stackrel{\mu}\rightarrow x\mu)(g))$. 
Since the coproduct functor is exact in  $\Ab$, the functor $\Lan^S$ 
is exact. Moreover, $Lan^S$ is left adjoint to the exact functor
$(-)\circ S: \Ab^{M^{op}} \rightarrow \Ab^{(M/X)^{op}}$ and 
therefore it carries projective objects to projective.
Consequently, $\Lan^S$ carries the projective resolution
$0\leftarrow F \leftarrow P_0 \leftarrow P_1 \leftarrow \cdots$ 
of the object $F\in \Ab^{(M/X)^{op}}$ to a projective resolution 
of
$\Lan^S\in \Ab^{M^{op}}$.
Applying the functor $\coLim^{(M/X)^{op}}$ to the resolution 
of the object $F$  we obtain a chain complex which homology groups 
are isomorphic to $\coLim_n^{(M/X)^{op}}F$.
It follows from an isomorphism
 $\coLim^{M^{op}}\Lan^S \cong \coLim^{(M/X)^{op}}$ 
that this chain complex may be obtained 
by the application of the functor
 $\coLim^{M^{op}}$ 
to a projective resolution of the right 
 $M$--module $\Lan^S F$. Hence its the homology groups are isomorphic 
to $\coLim_n^{M^{op}}\Lan^S F$.
Consequently, the groups $\coLim_n^{(M/X)^{op}}F$ are isomorphic to 
the homology groups of  $M$ with coefficients in the 
right  $M$--module $\Lan^S F$.

Substitute $M=M(E,I)$. By Theorem \ref{main} we obtain 
the chain complex for the calculating the homology groups 
with coefficients in $\Lan^S F$

$$
0 \leftarrow \Lan^S F \stackrel{d_1}\leftarrow
\bigoplus\limits_{a_1\in E} \Lan^S F
\stackrel{d_2}\leftarrow
\bigoplus\limits_{{a_1<a_2}\atop{(a_1,a_2)\in I}} \Lan^S F
\leftarrow  \cdots 
$$
$$
~~~~~~~~~~~~~~~~~~\cdots ~~~\leftarrow  
\bigoplus\limits_{{a_1<a_2<\cdots<a_{n-1}}\atop{(a_i,a_j)\in I} } 
\Lan^S F
\stackrel{d_n}\longleftarrow 
\bigoplus\limits_{{a_1<a_2<\cdots<a_n}\atop{(a_i,a_j)\in I} } 
\Lan^S F \leftarrow \cdots~,
$$
with differentials 
$$
d_n((a_1, \cdots, a_n),g) = 
\sum\limits_{s=1}^n (-1)^s ((a_1, \cdots, \widehat{a_s} , \cdots, a_n),
g\cdot a_s-g).
$$

We have $\Lan^S F=\bigoplus\limits_{x\in X}F(x)$ and the actions 
of $d_n$ are defined by values at  $(x,a_1\cdots a_n,g)$, $g\in F(x)$.
Hence this chain complex is isomorphic to the chain complex  
which members and differentials are described 
in the theorem.
\hfill $\Box$

\medskip
\noindent
{\bf Examples of calculating the homology groups 
of $M(E,I)$--sets.}
\begin{example}
Let $M(E,I)$ be a locally bounded free 
partially commutative monoid
which acts on the one-point set 
 $X=\{\star\}$.
In this case by Theorem 
 \ref{main2} the homology groups 
$\coLim_n^{(M(E,I)/X)^{op}}\Delta\ZZ$ are isomorphic 
to the homology of the chain complex
$$
0 \leftarrow 
 \ZZ \stackrel{d_1}\longleftarrow 
 \bigoplus\limits_{a_1}\ZZ \stackrel{d_2}\longleftarrow 
 \bigoplus\limits_{ 
{a_1<a_2}\atop{(a_1,a_2)\in I} } \ZZ
\leftarrow  \cdots 
$$
$$
~~~~~~~\cdots ~~~\leftarrow  
\bigoplus\limits_{ 
{a_1<a_2<\cdots<a_{n-1}}\atop{(a_i,a_j)\in I} } 
\ZZ
\stackrel{d_n}\longleftarrow 
\bigoplus\limits_{{a_1<a_2<\cdots<a_n}\atop{(a_i,a_j)\in I} } 
\ZZ \leftarrow \cdots~,
$$
where  
$$
d_n(a_1\cdots a_n) = 0
$$ 
Consequently, $\coLim_n^{(M(E,I)/\{\star\})^{op}}\Delta\ZZ\cong \ZZ^{\,p_n}$, 
where $p_n$ is the number of subsets
 $\{a_1, \cdots, a_n\} \subseteq E$ consisting of mutually 
commuting generators. (The number of empty subsets $p_0= 1$.)
\end{example}

We will show that the homology groups of $M(E,I)$--sets can have 
a torsion. 

Recall that a 
{\em simplicial schema} is a pair $(X,\mathfrak{M})$ consisting of
a set  $X$ with a set  $\mathfrak{M}$ 
of its nonempty subsets for which the following 
implication  holds 
$$
	S \subseteq T \in \mathfrak{M} \Rightarrow S \in \mathfrak{M}.
$$
Let $(X,\mathfrak{M})$ be a simplicial schema.
Consider an arbitrary total order relation 
 $<$ on $X$.
Let 
$$
X_n = \{ (x_0, x_1, \cdots, x_n) : \{x_0, x_1, \cdots, x_n\} 
\in \mathfrak{M}
\quad \& \quad x_0< x_1< \cdots< x_n \}.
$$
Consider the family of abelian groups
$$
C_n(X,\mathfrak{M})=\left\{
\begin{array}{ll}
L(X_n), & $for$ \quad  n \geq 0,\\
0, & $otherwise$.
\end{array}
\right.
$$
Define homomorphisms
$d_n: C_n(X,\mathfrak{M}) \rightarrow C_{n-1}(X,\mathfrak{M}) $
at elements of basis by
$$
d_n(x_0, x_1, \cdots, x_n) = 
\sum_{i=0}^n (-1)^i (x_0, \cdots, x_{i-1}, x_{i+1}, \cdots, x_n).
$$

It is well known that  $(C_n(X,\mathfrak{M}), d_n)$ 
is a chain complex.
Denote  $C_*(X,\mathfrak{M})=(C_n(X,\mathfrak{M}), d_n)$.
It easy to prove that homology groups of the chain complex 
 $C_*(X,\mathfrak{M})$ do not depend on the order relation. 
These groups are called {\em homology groups
$H_n(X,\mathfrak{M})$
of the simplicdial schema  $(X,\mathfrak{M})$}.

\begin{example}\label{exsimpsh}
Let $M(E,I)$ be any locally bounded free partially commutative monoid 
with the right action on the set $X=\{x_0, \star\}$ by  the lows 
$$
x_0\cdot a = \star ,\\
\star\cdot~ a = \star ~\mbox{для всех} ~ a\in E\,.
$$
Let $\ZZ[x_0]: (M(E,I)/X)^{op} \rightarrow \Ab$ be a functor 
such that
 $\ZZ[x_0](x_0)=\ZZ$ and $\ZZ[x_0](\star)=0$.
By Theorem \ref{main2} the groups
$\coLim_n^{(M(E,I)/X)^{op}}\ZZ[x_0]$ are isomorphic 
to the homology groups of the chain complex
$$
0 \leftarrow 
 \ZZ \stackrel{d_1}\longleftarrow 
 \bigoplus\limits_{a_1}\ZZ \stackrel{d_2}\longleftarrow 
 \bigoplus\limits_{ 
{a_1<a_2}\atop{(a_1,a_2)\in I} } \ZZ
\leftarrow  \cdots 
$$
$$
~~~~~~~\cdots ~~~\leftarrow  
\bigoplus\limits_{ 
{a_1<a_2<\cdots<a_{n-1}}\atop{(a_i,a_j)\in I} } 
\ZZ
\stackrel{d_n}\longleftarrow 
\bigoplus\limits_{{a_1<a_2<\cdots<a_n}\atop{(a_i,a_j)\in I} } 
\ZZ \leftarrow \cdots~,
$$
where 
$$
d_n(a_1\cdots a_n) = 
\sum_{s=1}^n(-1)^{s+1} 
(a_1\cdots \widehat{a_s}\cdots a_n)
$$ 
Consequently,
$\coLim_n^{M(E,I)/X^{op}} \ZZ[x_0] \cong H_{n-1}(E,\mathfrak{M})$
for all $n\geq 2$. 
\end{example}

Let $(M(E,I),X)$ be a pair consisting of 
a free commutative monoid and a right 
$M(E,I)$--set $X$.
Consider sets
$$
Q_n X = \{(x,a_1\cdots a_n): a_1<\cdots< a_n \&
(1\leq i < j \leq n \Rightarrow (a_i,a_j)\in I) \}
$$
 and maps
$$
 Q_n X {{{\partial^{n,0}_i}\atop{\longrightarrow}}\atop
	{{\longrightarrow}\atop{\partial^{n,1}_i}}}  Q_{n-1}X~,  
\quad n\geq 1~, ~~1\leq i\leq n~, ~~\varepsilon\in\{0,1\}~,
$$
defined by 
$$
{\partial^{n,0}_i}(x,a_1\cdots a_n) =
(x, a_1\cdots \widehat{a_i} \cdots a_n)~,
\quad
{\partial^{n,1}_i}(x,a_1\cdots a_n) =
(x a_i, a_1\cdots \widehat{a_i} \cdots a_n)~.
$$ 
\begin{lemma}
The sets $Q_n X$ and the maps
${\partial^{n,0}_i}$, ${\partial^{n,1}_i}$
constitute a precubical set.
\end{lemma}
{\sc Proof.} It may be directly checked that for every $x\in X$ 
and $\alpha, \beta \in \{0,1\}$ there are the following equalities 
$\partial_i^{n-1,\alpha} \partial_j^{n,\beta}(x) = 
\partial_{j-1}^{n-1, \beta}\partial_i^{n,\alpha}(x)$, 
for $1\leq i\leq n$. Hence 
$\partial_{i}^{n-1,\alpha} \partial_j^{n,\beta} = 
\partial_{j-1}^{n-1,\beta}\partial_i^{\alpha,n}$.
\hfill $\Box$

\begin{corollary}
Let  $X: M(E,I)^{op}\rightarrow \Set$ be a right $M(E,I)$--set
 Then the homology groups
$\coLim_n^{(M(E,I)/X)^{op}}\Delta\ZZ$ are isomorphic to integral 
homology groups of the precubical set
 $Q_* X$.
\end{corollary}
{\sc Proof.} Substitute in Theorem \ref{main2}  
$F=\Delta\ZZ$.
By Definition \ref{homolcub}
 the homology groups of the resulting chain complex  
 $C_*(Q_* X)$ 
are precisely the homology groups of the precubical set $Q_* X$.
\hfill $\Box$

\subsection{Homology of asynchronous transition systems}

\begin{definition}\label{defast}{\rm \cite{win1995}}
An {\em asynchronous transition system}
$$
A=(S,s_0,E,I,Tran)
$$
is a data consisting of sets $S$ and $E$,
a member $s_0\in S$, a subset $Tran\subseteq S\times E\times S$, 
and an irreflexive symmetric relation $I\subseteq E\times E$ 
satisfying to the following three conditions

(i) for each $e\in E$ there are  $s,s'\in S$
such that $(s,e,s')\in Tran$;

(ii) if $(s,e,s')\in Tran$ and $(s,e,s'')\in Tran$, then $s'=s''$;

(iii) for every pair $(e_1,e_2)\in I$ and for arbitrary triples 
 $(s,e_1,s_1)\in Tran$,
$(s_1,e_2,u)\in Tran$ there exists $s_2\in S$ such that
$(s,e_2,s_2)\in Tran$ and $(s_2,e_1,u)\in Tran$.

Members of $S$ are called {\em states}, members of $Tran$ are 
{\em transitions}, $s_0\in S$ is {\em an initial state}, $I$ is
{\em independence relation}.
\end{definition}

A {\em morphism of asynchronous transition systems}
$$
(\eta,\sigma): (S,s_0,E,I,Tran)\rightarrow (S',s'_0,E',I',Tran')
$$
consists of a map  $\sigma: S\rightarrow S'$
and a partial function $\eta: E\rightharpoonup E'$
for which $\sigma(s_0)=s'_0$ and the following two conditions hold

(i)  $(s,a,t)\in Tran \Rightarrow (\sigma(s),\eta(a), \sigma(t))\in Tran'$,
if $\eta(a)$ is defined and otherwise $\sigma(s)=\sigma(t)$;

(ii) if $(e_1,e_2)\in I$ and  $\eta(e_1)$, $\eta(e_2)$
are both defined, then
$(\eta(e_1),\eta(e_2))\in I'$.

The category of asynchronous transition systems we denote by
$\fA$.

\medskip
\noindent
{\bf Asynchronous transition systems and right pointed sets over monoids.}

\begin{definition}
A {\em right pointed set over a monoid} is a triple
$(M,\cdot,X)$ consisting of a monoid  $M$, a right pointed $M$--set
$X$, and a map
$\cdot: X\times M \rightarrow X$ such that

${\star}\cdot\mu={\star}$ for all $\mu\in M$;

$(x\cdot \mu_1)\cdot\mu_2= x\cdot(\mu_1 \mu_2)$ for all $x\in X$,
$\mu_1, \mu_2\in M$;

$x\cdot 1=x$, for all $x\in X$.

\end{definition}
Here $1$ is the identity element of the monoid.
The symbol $'\cdot\,'$ is omitted usually.
So a right pointed set over a monoid is denoted by  $(M, X)$.

A {\em morphism} of right pointed sets over monoids 
$(\eta,\sigma): (M,X) \rightarrow (M',X')$ is a pair consisting of 
a monoid homomorphism  $\eta:M\rightarrow M'$ and based map
 $\sigma: X\rightarrow X'$ satisfying 
$\sigma(x\cdot \mu)= \sigma(x)\cdot\eta(\mu)$ for all $x\in X$ and $\mu\in M$.
Denote the category of right pointed sets over monoids by $(Mon, \Set_*)$.

Define a map $Ob(\fA) \rightarrow (Mon, \Set_*)$ assigning to every 
asynchronous transition system 
$(S,s_0,E,I,Tran)$ a pair $(M(E,I), S_*)$ where 
$S_*$ is the pointed set 
with the action 
$(-)\cdot{e}: S\cup \{\star\}\rightarrow S\cup \{\star\}$ 
on generators $e\in E$ by
$$
s\cdot e=\left\{
\begin{array}{cl}
s'~, & \mbox{if}   ~~  (s,e,s')\in Tran , \\
\star~, & \mbox{if there is not} ~s'~ \mbox{such that} 
 ~~  (s,e,s')\in Tran \,.
\end{array}
\right.
$$

 It is proved in \cite{X20042} that 
the map
$Ob(\fA) \rightarrow (Mon, \Set_*)$ may be extended  
 to a functor $\fA \rightarrow (Mon, \Set_*)$.
This functor assigns to any morphism 
$$
(\eta,\sigma): (S,s_0,E,I,Tran)\rightarrow (S',s'_0,E',I',Tran')
$$
of asynchronous transition systems the morphism 
$(\eta^*,\sigma_*): (M(E,I),S_*) \rightarrow (M(E',I'),S'_*)$
consisting of an extension 
of the map $\eta: E\rightarrow E'$ to a monoid homomorphism 
and of natural transformation 
$\sigma_*: S_*\rightarrow S'_*\circ (\eta^*)^{op}$
consisting of
the map 
$$
\sigma_*(s)=\left\{
\begin{array}{cl}
\sigma_(s)~, & \mbox{if } ~s\in S,\\
\star~, & \mbox{if } s=\star.
\end{array}
\right.
$$

A {\em forgetful} functor $U: \Set_*\rightarrow \Set$ defined 
by $X\mapsto X$ and $f\mapsto f$. 
An {\em augmented state category} 
$K_*(T)$ \cite[Definition 5.1]{X20042} 
of an asynchronous transition system $T=(S,s_0,E,I,Tran)$
is the category which  objects are the members 
$s\in S\cup \{\star\}$ and the morphisms 
are triples $(\mu, s, s')$ consisting of 
 $\mu\in M(E,I)$, $s,s'\in S$ satisfying to $s\cdot\mu=s'$.

A {homology groups $H_n(K_*(T),F)$ of an asynchronous transition 
system $T$ with coefficients in a functor
 $F:K_*(T)\rightarrow \Ab$} are abelian groups
 $\coLim_n^{K_*(T)}F$.

Events $e_j\in E$, $j\in J$ are said to be {\em mutually independent},
if  $(e_j,e_{j'})\in I$ for all  $j,j'\in J$ such that $j\not=j'$.

\begin{theorem}\label{main3}
Let $T=(S,s_0,E,I,Tran)$ be an asynchronous transition system
without infinite sets of mutually independent 
events with an arbitrary total order relation $'<\,'$ on $E$ and   
$F: K_*(T) \rightarrow \Ab$ a functor. Then 
the homology groups $H_n(K_*(T),F)$ are isomorphic to homology groups
of the following chain complex
$$
0 \leftarrow 
\bigoplus\limits_{s\in S\cup \{\star\}} F(s) \stackrel{d_1}\longleftarrow 
\bigoplus\limits_{(s,e_1)} F(s)
\stackrel{d_2}\longleftarrow \bigoplus\limits_{{(s, e_1 e_2), e_1<e_2}
\atop{(e_1,e_2)\in I}} F(s)
\leftarrow  \cdots 
$$
$$
\cdots ~~
\leftarrow  
\bigoplus\limits_{(s, e_1\cdots e_{n-1}), 
{e_1<e_2<\cdots<e_{n-1}}\atop{(e_i,e_j)\in I} } 
F(s)
\stackrel{d_n}\longleftarrow 
\bigoplus\limits_{{(s, e_1\cdots e_n),
e_1<e_2<\cdots<e_n}\atop{(e_i,e_j)\in I} }
F(s) \leftarrow \cdots~ 
$$
where 
$$
d_n(s,e_1\cdots e_n,g) = 
$$
$$
\sum_{i=1}^n(-1)^i \left(
(s\cdot e_i, e_1\cdots \widehat{e_i}\cdots e_n, 
F(s\stackrel{e_i}\rightarrow s\cdot e_i)g)-
(s, e_1\cdots \widehat{e_i}\cdots e_n, g) \right)
$$ 
\end{theorem}
{\sc Proof.} It easy to see that $K_*(T)$ is isomorphic to 
a category which is dual to $M(E,I)/U\circ S_*$.
Therefore the wanted assertion follows from Theorem \ref{main2}.
\hfill $\Box$

If the asynchronous transition system
 $T$ has  $\leq n$ mutually independent events, 
then the chain complex obtained in Theorem 
 \ref{main3}, 
has the length $\leq n$. Hence we get the following 
\begin{corollary}
Let $n>0$ be the maximal number of mutually 
independent events in an asynchronous 
transition system  $T$. 
Then $H_k(K_*(T),F)=0$ for all $k>n$.
\end{corollary}
This assertion gives the positive answer to  
my question  
\cite[Open problem 2]{X20042}.  
Theorem \ref{main3} leads to the algorithm 
of calculating the integral homology groups of an asynchronous 
transition system with  
finite sets  
$S$ and $E$ by the computations of the Smith normal 
form of the differentials $d_n$.

\end{document}